\documentclass[11pt]{amsart}

\usepackage{epigamath}

\usepackage[english]{babel}

\numberwithin{equation}{section}

\usepackage[shortlabels,inline]{enumitem}

\usepackage{euscript,amsmath,amssymb,mathabx}
\usepackage[arrow,matrix,curve]{xy}
\usepackage{array,longtable,graphicx,enumitem,url,multirow}

\newtheorem{m-theorem}{Theorem}[section]
\newtheorem{m-corollary}[m-theorem]{Corollary}
\newtheorem{m-lemma}[m-theorem]{Lemma}
\newtheorem{m-proposition}[m-theorem]{Proposition}

\newtheorem*{thm-nono}{Theorem}
\newtheorem*{prop-nono}{Proposition}
\newtheorem*{lm-nono}{Lemma}
\newtheorem*{conj-nono}{Conjecture}

\theoremstyle{definition}

\newtheorem{m-definition}[m-theorem]{Definition}
\newtheorem{m-notation}[m-theorem]{Notation}
\newtheorem{m-asmp}[m-theorem]{Assumptions}

\newtheorem*{def-nono}{Definition}

\theoremstyle{remark}

\newtheorem{m-example}[m-theorem]{Example}
\newtheorem{m-remark}[m-theorem]{Remark}
\newtheorem{m-question}[m-theorem]{Question}

\newenvironment{m-eqn}
  {\begin{equation}}
  {\end{equation}}

\setlist[enumerate,1]{label={\rm(\roman*)}, ref={\rm\roman*}}
\setlist[enumerate,2]{label={\rm(\alph*)}, ref=\rm\alph*} 

\newlist{aenum}{enumerate}{2}
\setlist[aenum,1]{label={\rm(\alph*)}, ref={\rm\alph*}}
\setlist[aenum,2]{label={\rm(\roman*)},ref={\rm\roman*}}

\makeatletter\newcommand\customitem[1]{\item[#1{\MakeLinkTarget[item]{}}]\def\@currentlabel{#1}}\makeatother

\let\euf\EuScript 
\let\cal\mathcal
\let\mbb\mathbb

\newcommand{\Ad}{{\rm Ad}}
\newcommand{\BD}{{O_1}}
\let\dta\delta

\let\ges\geqslant 
\let\les\leqslant 
\let\hra\hookrightarrow 

\let\nit\noindent
\let\ovl\overline 
\let\sm\setminus 
\let\Si\Sigma 
\let\srel\stackrel 
\let\tld\tilde 
\let\unbar\underbar 
\let\vpi\varpi 
\let\vphi\varphi 
\let\disp\displaystyle 

\newcommand{\cd}{\mathop{\rm cd}\nolimits}
\newcommand{\codim}{\mathop{\rm codim}\nolimits}
\renewcommand{\dim}{\mathop{\rm dim}\nolimits}
\newcommand{\etl}{{\rm\acute{e}t}}
\newcommand{\gen}{{\rm gen}}
\newcommand{\iris}{{\rm iris}}
\newcommand{\lcit}{{\textit{loc.\,cit.}}}
\newcommand{\ocit}{{\textit{op.\,cit.}}}
\newcommand{\kk}{{\mbb C}}

\newcommand{\eI}{{\euf I}}

\newcommand{\eL}{{\euf L}}
\newcommand{\pos}{{>0}}

\newcommand{\rd}{{\rm d}}
\let\rst\vert 
\newcommand{\eN}{{\euf N}}
\newcommand{\eO}{{\euf O}}
\newcommand{\eT}{{\euf T}}
\newcommand{\cY}{{\cal Y}}

\newcommand{\VV}{V}
\newcommand{\xx}{x}
\newcommand{\XX}{X}
\newcommand{\YY}{Y}
\newcommand{\tpl}{{\rm top}}
\newcommand\Hom{\mathop{\rm Hom}\nolimits}
\newcommand\Pic{\mathop{\rm Pic}\nolimits}
\newcommand{\Spec}{\mathop{\rm Spec}\nolimits}
\newcommand{\Sym}{\mathop{\rm Sym}\nolimits}
\newcommand{\res}{\mathrm{res}}
\newcommand{\proj}{\mathrm{proj}}
\newcommand{\open}{\mathrm{open}}
\newcommand{\num}{\mathrm{num}}
\newcommand{\kod}{\mathrm{kod}}
\newcommand{\rmdiv}{\mathrm{div}}
\newcommand{\amp}{\mathrm{amp}}
\DeclareMathOperator{\Grass}{Grass}
\DeclareMathOperator{\Eff}{Eff}
\DeclareMathOperator{\rank}{rank}
\DeclareMathOperator{\Sing}{Sing}

\newcommand{\lra}{\longrightarrow}

\EpigaVolumeYear{10}{2026} \EpigaArticleNr{18} \ReceivedOn{January 13, 2024}
\InFinalFormOn{February 26, 2026}
\AcceptedOn{April 13, 2026}

\title{G3-criteria and applications}
\titlemark{G3-criteria and applications}

\author{Mihai Halic}
\address{EMR-matrix, Arad, Romania}
\email{mihai.halic@gmail.com}

\authormark{M.~Halic}

\AbstractInEnglish{The G3-property of a subvariety was introduced by Hironaka--Matsumura and plays an important role for deducing connectedness and extension results. Unfortunately, it is a rather elusive notion, which is not always easy to establish. Most of the existing work is concentrated on subvarieties of homogeneous varieties.

  The first goal of this article is to show that mobility assumptions on the subvariety, considered in works of B\u{a}descu, Chow, Debarre, Voisin, yield a certain partial positivity property, slightly stronger than G3, previously introduced by the author. Second, we apply the result to prove that, in numerous situations, the splitting of the normal bundle of a smooth two-codimensional subvariety implies that it is a complete intersection.}



\begin{document}



\maketitle

\begin{prelims}

\DisplayAbstractInEnglish

\bigskip

\DisplayKeyWords

\medskip

\DisplayMSCclass

\end{prelims}


\newpage

\setcounter{tocdepth}{1}

\tableofcontents


\maketitle

\section{Introduction}

Hironaka--Matsumura~\cite{hir-mat} introduced generating properties (G1, G2, and G3) for a subvariety $\YY$ of a projective variety $\XX$ in order to quantify the amount of geometric information carried by a formal (or tubular) neighbourhood of the subvariety. The G3-property is the strongest: the formal completion along~$\YY$ birationally determines $\XX$. The G2-property requires only that the ambient space is determined up to a finite cover. 

The author~\cite{hlc-g2} proved that a local complete intersection subvariety with non-pseudo-effective co-normal bundle possesses the G2-property; this generalizes the classical result of Hartshorne~\cite{hart-cd}, which assumes that the normal bundle is ample. Moreover, the author showed that the non-pseudo-effectiveness property holds for strongly movable subvarieties. Combined, these facts yield a practical geometric G2-criterion.

In contrast, the G3-property is rather elusive, in spite of being important for addressing connectedness and extension issues. Most of the existing literature is focused on studying the properties of (arbitrary) subvarieties of homogeneous spaces (rational or abelian); see~\cite{bad,chow,dbar-mumford,falt-homog}. There are few sufficient geometric criteria for it; one of the best known is due to Speiser; see \cite{hart-as}. The author~\cite{hlc-pas} introduced the notion of a partially ample subvariety and proved that the G3-property holds for subvarieties satisfying the weakest such partial amplitude condition (named $1^\pos$, one-positive); the interest of the partial amplitude property was shown by proving Fulton--Hansen-type connectedness results for pre-images. The $1^\pos$-property has, however, an unpractical aspect: it involves the cohomological dimension of the complement $\XX\sm\YY$, which is not straightforward to control (see~\cite{hlc+taj} for several methods to estimate it).

The main goal of this article is to give manageable sufficient, geometric conditions for the $G3$-property of lci and smooth subvarieties --~this suffices for numerous situations~-- and to present a variety of applications. The key notion used  here is the mobility of $\YY$ in $\XX$; it allows one to deduce the G3- from the G2-property. The mobility feature is already present in works of B\u{a}descu~\cite{bad}, Chow~\cite{chow}, Debarre~\cite{dbar-mumford}, Voisin~\cite{voisin-coniv}; the connection between mobility (of curves) and pseudo-effectiveness (of line bundles) is the main topic of~\cite{bdpp}. 

The article is structured as follows: we first recall relevant background and then analyse several mobility conditions for families $\cY$ of subvarieties $\YY\subset\XX$. In analogy to the case of group actions, the fibres of $\cY\to\XX$ are called isotropy spaces; they parametrize pointed deformations within the family. It turns out that the irreducibility of the general isotropy space plays an important role;  this property is automatic for almost homogeneous varieties, and it is frequently satisfied for simply connected $\XX$ (\textit{cf.}~Proposition~\ref{prop:irr-fibre}).

The first main result, Theorem~\ref{thm:g3}, generalizes a classical result of Chow~\cite{chow}: briefly, if $\YY$ intersects all effective divisors of $\XX$ and moves sufficiently within a family with irreducible general isotropy, then it has the G3-property. Let us mention two applications. 
The latter is a positive answer to a question raised by Debarre~\cite[Conjecture 2.9]{dbar-mumford}.

\begin{thm-nono}{}
Suppose $\YY$ is lci, intersects numerically non-trivially all effective divisors in $\XX$, and is strongly movable in a family $\cY$ whose general isotropy is irreducible. Then the following statements hold true:
\begin{enumerate}
\item {\rm(\textit{cf.}~Corollary~\ref{cor:mori})}
  The subvariety $\YY\subset\XX$ is $1^\pos$, hence G3, in either one of the following situations:
  
\begin{enumerate*}
\item $\XX$ is rationally connected; or\, 
\item $\YY\!$ has ample normal bundle and $2\dim\YY\!\ges\dim\XX$.
\end{enumerate*}
\item {\rm(\textit{cf.}~Corollary~\ref{cor:dbar})}
There exist an \'etale cover $\XX'\srel{f}{\to}\XX$ and a subvariety $\YY'\subset\XX'$ such that $\YY'\srel{f}{\to}\YY$ is an isomorphism and $\YY'$ is $G3$ in $\XX'$.
\end{enumerate}
\end{thm-nono}

Further applications concern two-codimensional subvarieties. Recall that, conjecturally, the only two-codimensional subvarieties of the projective space $\mbb P^n$, $n\ges6$, are complete intersections. Faltings~\cite{falt-krit} proved that if the normal bundle of $\YY\subset\mbb P^n$ splits and $\dim\YY\ges n/2$, then $\YY$ is indeed a complete intersection. 
The starting point of our investigation is that one may ask the analogous question for \emph{arbitrary} ambient varieties; we consider the case when $\YY\subset\XX$ has split normal bundle. The main result is formulated in Theorem~\ref{thm:g3-ci}; we state here a brief application concerning subvarieties of rational homogeneous spaces.

\begin{thm-nono}[\textit{cf.}~Corollary~\ref{cor:rtl}]
Let $G$ be a simple linear algebraic group with $\rank(G)\ges6$, $P\subset G$ a maximal parabolic subgroup, and $\XX=G/P$. Then any two-codimensional, smooth subvariety of\, $\XX$ with split normal bundle is complete intersection. A similar result holds for $\YY\!$ lci $($cf.\ Theorem~\ref{thm:ci-sing}\,$)$.
\end{thm-nono}

Although there is a rich literature devoted to the complete intersection property for subvarieties of projective spaces, to the author's knowledge, there are basically no references investigating the analogous issue for subvarieties of homogeneous varieties. 

\subsection*{Acknowledgments}
I thank the referee(s) for the careful reading and comments which helped improving the presentation, for pointing out  inaccuracies in the original manuscript which motivated the considerations in Section~\ref{ssct:irr}, and for the question on the lci case leading to Section~\ref{ssct:lci}. 


\section{Movable and generating subvarieties}\label{sct:gen-fam}

\subsection{Background}\label{ssct:setup}

In this article, $\XX$ stands for a complex, smooth projective variety. For a connected subvariety $\YY$ defined by the ideal $\eI_\YY\subset\eO_\XX$, let $\hat \XX_\YY$ be the formal completion of $\XX$ along $\YY$; its structure sheaf is ${\varprojlim}_m\,\eO_\XX/\eI_\YY^{m}$. The rational formal functions $\kk(\hat \XX_\YY)$ on $\XX$ along $\YY$ form a field. 
Hironaka--Matsumura~\cite{hir-mat} introduced the following generating properties: 
\begin{itemize}
\item $\YY$ is G2 in $\XX$ if $\kk(\XX)\hra \kk(\hat \XX_\YY)$ is a finite field extension. 
\item $\YY$ is G3 in $\XX$ if $\kk(\XX)\hra \kk(\hat \XX_\YY)$ is an isomorphism.
\end{itemize}
For a detailed treatment of these topics, the reader may consult the reference B\u{a}descu~\cite{bad}. The following well-known result due to Speiser relates the G2- and G3-properties.

\begin{thm-nono}[\textit{cf.}~\cite{hart-as}]
The following conditions are equivalent $($for a closed subscheme $\YY\subset\XX)$: 
\begin{enumerate}
\item $\YY$ is G3 and $\YY$ intersects all divisors of $\XX$;
\item $\YY$ is G2 and $\cd(\XX\sm \YY)\les\dim\XX-2$. 
$($Here $\cd(\cdot)$ stands for ``cohomological dimension''.$)$
\end{enumerate}
\end{thm-nono}

The author~\cite{hlc-pas} introduced a sequence of partial amplitude conditions for a subvariety, building upon work of Ottem~\cite{ottm}; the weakest was called $1^\pos$ (one-positivity). 

\begin{m-definition}[\textit{cf.}~{\cite[Proposition 1.6(iii)]{hlc-pas}}]
A locally complete intersection (lci) subvariety $\YY$ is $1^\pos$ if $\,\eN_{\YY/\XX}^\vee\!=\eI_\YY/\eI_\YY^2\;\text{is not pseudo-effective and}\;\cd(\XX\sm\YY)\les\dim\XX\!-2.$
\end{m-definition}

Recall, see~\cite[Corollary 2.6]{hlc-g2}, that the non-pseudo-effectiveness of $\eN_{\YY/\XX}^\vee$ implies that $\YY$ is $G2$ in~$\XX$. Together with Speiser's criterion, this implies that lci, $1^\pos$ subvarieties possess the G3-property. The drawback of this definition is that it requires controlling the cohomological dimension of $\XX\sm\YY$, which is not straightforward. (See~\cite{hlc+taj} for various methods to estimate it.)

To motivate the forthcoming mobility considerations, let us also mention that the author proved, see~\cite[Proposition 3.5]{hlc-g2}, that strongly movable lci subvarieties --~this is a simple geometric condition, see below~-- have non-pseudo-effective co-normal bundle. 


\subsection{Mobility}

For curves and divisors, the concept is classical. In higher codimensions, the notion of (very) moving subvarieties was introduced by Voisin~\cite{voisin-coniv}, in the attempt to geometrically characterize big subvarieties. A notion with similar flavour, that of generating subvarieties, appears in Chow~\cite{chow}; see also Debarre~\cite{dbar-mumford}, B\u{a}descu~\cite{bad}. The latter notion turns out to be essential for investigating the G3-property of subvarieties of homogeneous varieties. With this motivation, we consider the following. 

\begin{m-definition}\label{def:YY-gener}\leavevmode
\begin{aenum}
\item\label{d:YY-g-a} Let $S$, the parameter space, be a quasi-projective variety, and let 
\begin{equation}
\xymatrix@R=0.49em@C=3em{
&\cY\ar[dl]_-\vpi\ar[dr]^\rho &&\cY\subset S\times\XX\\ S&&\XX&}\label{eq:syx}
\end{equation}
be an irreducible, flat $S$-family of subvarieties of $\XX$; we assume that $\cY$ is normal. We consider the following varieties and schemes (they are of interest whenever non-empty):
\begin{enumerate}
\item\label{d:YY-g-a1} Let $\Si_1:=S$, $\cY_{1}=\cY$, $\vpi_1:=\vpi$, and let $\Delta_\cY\subset\cY\times\cY$ be the diagonal.
  
\item\label{d:YY-g-a2} The incidence variety $\Si_2\subset S\times S$ is the closure of the locus of the following pairs $(s,s')$: first, one has $\YY_s\cap\YY_{s'}\neq\emptyset$, $\YY_s$, $\YY_{s'}$; second, $(s,s')$ belongs to an irreducible component containing $(s,s), (s',s')$; that is,  $\YY_s\cap\YY_{s'}$ deforms (non-equidimensionally) to the ``identity intersections'' $\,\YY_s\cap\YY_s$, $\YY_{s'}\cap\YY_{s'}\,$. Below is the formal definition. 

Let $\mbb S:=\{ (s,s')\mid \YY_s=\YY_{s'} \}$; if $S$ parametrizes distinct subvarieties of $\XX$, then $\mbb S$ is the diagonal  $\Delta_S\subset S\times S$. We define $\Sigma_2$ to be the union of those irreducible components
\[
  \xymatrix{R\subset\overline{(\vpi_1,\vpi)(\cY_1\times_\XX\cY)\;\sm\;\mbb S}\subset S\times S\ar@<.35ex>[r]^-{\proj_1}\ar@<-.35ex>[r]_-{\proj_2}&S}\quad
\text{such that }
\proj_1(R)=\proj_2(R)\subset S
\](we denote the common image by $C$) and $R\supset \Delta_C$ (the diagonally embedded $C$). 
The reason for eliminating $\mbb S$ is to consider only those $\YY_s$, $s\in S$, which possess non-trivial, intersecting deformations. For movable families (see below), this definition simplifies considerably (see Remark~\ref{rmk:G}(\ref{rmk:G-2})).  

Let $\cY_2$ be defined as in the diagram below: 
\vskip-.25ex
\begin{equation}
\xymatrix@C=3.49em@R=1.79em{
\cY_2:=
\Si_2\times_{S\times S}(\cY_1\times\cY)  \ar[r]\ar@/^3ex/[rrr]|{\;\rho_2\;}\ar[d]\ar[dr]|{\;(\vpi_1,\vpi_2)\;}&
\cY_1\times\cY\ar[r]_<(.35){(\rho_1,\rho)}\ar[d]^<(.45){(\vpi_1,\vpi)}& \XX\times\XX\ar[r]_-{\proj_2\;}&\XX.
\\ 
\Si_2\ar[r]^<(.25){\iota}&S\times S&
}\label{eq:cY}
\end{equation} 
For $o\in S$, let $\{o\}\times \Sigma_o:=\iota^{-1}(\{o\}\times S)$ and $\cY_o:=\cY\rst_{\Sigma_o}$, and let $\rho_o:=\rho\rst_{\cY_o}$ --~the symbol ``$\rst$'' stands for restriction~-- be the evaluation morphism $\rho$ on those $\YY_s$ which intersect $Y_o$ non-trivially.
 
\item\label{d:YY-g-a3} Suppose $\Si_t\subset S^t$, $\cY_t=(\vpi_1,\dots,\vpi_t)^{-1}(\Si_t)\subset\cY^t$ are constructed, with morphisms 
\[
(\vpi_1,\dots,\vpi_t,\rho_1,\dots,\rho_t)\colon \cY_t\lra (S\times\dots\times S)\times (\XX\times\dots\times \XX).
\]
Here $\vpi_i$, $\rho_i$, $i=1,\dots, t$, are, respectively, the compositions of $\vpi$, $\rho$ with $\cY_t\to\cY^t\xrightarrow{\proj_i}\cY$. A~diagram similar to~\eqref{eq:cY} above allows us to define 
\begin{align*}
&\Si_{t+1}:=(\Sigma_t\times S)\times_{\Sigma_{t-1} \times S^{2}}\left(\Sigma_{t-1}\times\Sigma_2\right)\subset S^{t+1}, 
\\
&\cY_{t+1}:=\left(\vpi_1,\dots,\vpi_{t+1}\right)^{-1}\left(\Si_{t+1}\right)\subset\cY^{t+1}.
\end{align*}
There is a forgetful morphism of the $(t+1)^{\rm st}$ component $\Sigma_{t+1}\to\Sigma_{t}$, and also a ``repeating the last component'' morphism $\Sigma_t\to\Sigma_{t+1}$. We call $\Si_t$ the \emph{$t$-incidence variety of\, $S$}; roughly, it consists of sequences $(s_1,\dots,s_t)$ whose successive terms belong to $\Sigma_2$. A \emph{$\cY$-chain of length $t$} is a sequence $\YY_{s_1},\dots,\YY_{s_t}$, with $(s_1,\dots,s_t)\in\Si_t$. 
\end{enumerate}

\item \label{d:YY-g-b}
We define the following properties for the family $\cY$: 
\begin{enumerate}
\item\label{d:YY-g-b1} The family $\cY$ is \emph{movable} if $\rho\colon\cY\to \XX$ is dominant. 

\item\label{d:YY-g-b2} The family $\cY$ is \emph{strongly movable} if $(\vpi_1,\rho_2)\colon\cY_2\to S\times \XX$ is dominant and its smooth locus does intersect $(\varpi_1,\varpi_2)^{-1}(\Delta_S)$ above some open subset $S'$ of $S$. 

  In this situation, for $o\in S'$ and a general point $q\in\YY_o$, $\rd\rho_{o,q}\colon\eT_{\cY_o,(o,q)}\to\eT_{\XX,q}$ is surjective and $\rho_o$ is dominant. We say that $\YY_o$ is strongly movable. 

The family $\cY$ is \emph{strongly movable in codimension one} (codim-one, for short) if there is an open $O\subset\XX$  such that  
$O\subset\rho_2(\cY_2|_{\Sigma_2\cap (S'\times S)})$ and $\codim_\XX(\XX\sm O)\ges2$. 

\item\label{d:YY-g-b3} The family $\cY$ is \emph{generating} if $\big\{(s,x)\in S\times\XX\mid \exists t\ges 1,\;(s,x)\in(\vpi_1,\rho_t)(\cY_t)\big\}$ contains an open subset of $S\times\XX$. Equivalently, the last component of $\cY$-chains starting with $Y_o$, where $o\in S$ is general, sweeps out an open subset $O_o\subset\XX$.

\item\label{d:YY-g-b4} The family $\cY$ is \emph{chain connecting in codim-one} if the family $\cY$ is generating and the open subset above has the property $\;\codim_\XX(\XX\sm O_o)\ges2$.
\end{enumerate}

\item \label{d:YY-g-c}
We say that \emph{the smooth subvariety $Y\subset X$ is very moving} (\textit{cf.}~Voisin \cite[Definition 0.5]{voisin-coniv}) if the following condition holds:
\[
\begin{array}{l}
\forall\,[\text{generic}\;p\in\YY\;\text{and generic}\;T\subset\eT_{\XX,p}\;\text{of dimension}\;\dim\YY],\\[.5ex] 
\exists\,\text{smooth, nearby, embedded deformation $\YY'$ of $\YY$, with $p\in\YY'$, such that $\eT_{\YY',p}=T$.}
\end{array}
\]
The points $p$ constitute the \emph{very moving locus of $\YY$}.  Here we agree that the parameter space~$S$ parametrizes sufficiently many deformations $\YY'$ of $\YY$ such that the tangent spaces $\eT_{\YY', p}$, for $\YY'$ in $S_p:=\rho^{-1}(p)$, satisfy the property above. (Note that $S$ does not necessarily consist of all possible embedded deformations of $\YY$ in $\XX$. In Section~\ref{ssct:subvar-homog}, we consider only translates for a group action.) 
Voisin's condition means that 
\[
\,S_p\srel{\dta_p}{\lra} \Grass\left(\dim\YY,\eT_{\XX,p}\right),\quad \YY'\longmapsto \eT_{\YY',p} 
\]
is dominant about $\YY\in S_p$. Therefore, it remains so for $\YY'\in S$ in a neighbourhood of $\YY$; all small deformations $\YY'$ of $\YY$ are still very moving. By abuse of language, we say that \emph{the family $\cY$ is very moving} (after possibly shrinking $S$).

If $S_p$ is smooth at $\YY$, it is enough to require  $\eT_{S_p,\varpi(\YY)}\xrightarrow{\dta_{\varpi(\YY),p}}\Hom(\eT_{\YY,p},\eN_{\YY/\XX,p})$ to be surjective. Since $\eT_{S_p,\varpi(\YY)}\subset H^0(\YY,\eI_{\YY,p}\otimes\eN_{\YY/\XX})$, this map is simply the homomorphism
\[
H^0(\YY,\eI_{\YY,p}\otimes\eN_{\YY/\XX})\lra \Hom\left(\eT_{\YY,p},\eN_{\YY/\XX,p}\right),\quad \sigma\longmapsto\rd\sigma_p.
\]
Let us recall, see \cite[Remark 2.6]{voisin-coniv}, that one may assume, moreover, that $\eN_{\YY/\XX}$ is ample. This is still an open condition under deformations.

\item\label{d:YY-g-d}
The \emph{family $\cY$ is smooth at $p\in\XX$} if $\cY$ is smooth at every $(\varpi(\YY),p)\in\rho^{-1}(p)$ and $\rd\rho_{(\varpi(\YY),p)}$ is surjective: 
\[
\disp\sum_{\sigma\in\eT_{S_p,\varpi(\YY)}}\kern-2ex\rd\sigma_p\left(\eT_{\YY,p}\right)=\eN_{\YY/\XX,p}.
\]
We say that $p$ belongs to the \emph{smooth locus of\, $\cY$}. In complex-analytic language, this means that $p$ is a regular value of $\rho$. We say that $\rho$ is {\em smooth over codimension one points of\, $\XX$} --~for short, \emph{$\rho$ is smooth in codim-one}~-- if $\rho(\cY)$ contains an open subset $O$ such that $\codim_\XX(\XX\sm O)\ges2$ and $\cY$ is smooth at all $p\in O$. 

\item \label{d:YY-g-e}
  For short, we write $A^\gen$ for some open (general) subset of a variety or scheme $A$.
  \end{aenum}
\end{m-definition}

\begin{m-remark}\label{rmk:G} \leavevmode
  \begin{enumerate}
  \item\label{rmk:G-1} The formulation ``deformations of $\YY$'' refers \emph{only} to members of some family $\cY\to S$ containing $\YY$. The family is either given or suitably chosen, depending on the context.

The flatness of $\cY\to S$ is required whenever we invoke the deformation invariance of intersection products in $H^*(X)$ with members of the family. Otherwise, the definitions make sense in general. 

\item\label{rmk:G-2} We are interested in families $\cY$ which are movable, at the very least. The differential of $\rho$ is surjective ($\rho$ is smooth) on an open locus $\cY^\gen\subset\cY$; for any point $p\in\cY^\gen$, $\cY\times_\XX\cY$ is smooth about $(p,p)$. So there is only one irreducible component of $\cY^\gen\times_\XX\cY^\gen$ containing the diagonal. Its projection to both $S$-factors is dominant; let $\mbb\YY_2$ be its closure. Thus we get a ``natural''  component $\overline{(\vpi_1,\vpi_2)(\mbb\YY_2)\sm\mbb S}\subset\overline{(\vpi_1,\vpi_2)(\mbb\YY_2)\sm\Delta_S}\subset\Sigma_2\subset S^2$. If $Y_s$, $s\in S$, are pairwise-distinct subvarieties of $\XX$, then $\mbb S=\Delta_S$.

\item\label{rmk:G-3} Note that $\Sigma_2=S\times S$ as soon as the self-intersection satisfies $[Y]^2\neq0$; any two members of the family intersect non-trivially. Thus we have $\cY_t=\cY^t$ and $\Sigma_t=S^t$ for $t\ges1$.

  The self-intersection is non-trivial in this situation: for $o\in S^\gen$, $\YY_o$ is lci, $\eN_{\YY_o/\XX}$ is ample, and  $2\dim\YY_o\!\ges\dim\XX$. The excess intersection theorem, see~\cite{fu-laz}, yields $\deg_{\eO_\XX(1)}[Y]^2\!>0$. 

\item\label{rmk:G-4} Suppose $\XX$ is homogeneous for the action of a group $G$. Chow, see~\cite{chow, dbar-mumford, bad}, introduced the notion of a generating subvariety $\YY\subset\XX$. One verifies that his definition coincides with ours for $S:=G$ and $\cY:=G\times \YY$ (see Lemma~\ref{lm:strg-gen}).
  \end{enumerate}
\end{m-remark}

\begin{m-lemma}\label{lm:mov}\leavevmode
\begin{enumerate}
\item\label{l:mov-1} For the family $\cY$, we have the following implications: 
\[
\text{movable} \;\;\Longleftarrow\;\;\text{generating}\;\;\Longleftarrow\;\;\text{strongly movable}\;\;\Longleftarrow\;\;\text{very moving.}
\]
\item\label{iti:Ygen} 
The family $\cY$ generates $\XX$ if and only if\, $(\vpi_1,\rho_t)\colon\cY_t\to S\times\XX$ is dominant for some $t\ges 2$.
\item\label{l:mov-3}  Let $o\in S$ and $p_o,q_o\in\YY_o\subset\cY$ be smooth points.  
\begin{itemize}
\item If\, $\widehat{\rd\rho}_{p_o}\colon \eT_{\cY, p_o} \to \eN_{\YY_o/\XX,p_o}$ is surjective, then $\cY$ is movable. 
\item If\, $\widehat{\rd\rho}_{p_o}\oplus\widehat{\rd\rho}_{q_o}\colon  \eT_{\cY\times_S\cY, (p_o, q_o)} \to\eN_{\YY_o/\XX,p_o}\oplus\eN_{\YY_o/\XX,q_o}$ 
is surjective, then $\cY$ is strongly movable. $($Note that $\eT_{\cY\times_S\cY, (p_o, q_o)}= \{(\xi,\xi')\mid\rd\vpi_{p_o}(\xi)-\rd\vpi_{q_o}(\xi')=0$.$)$
\end{itemize}
Here $\widehat{\rd\rho}$ is the composition of $\rd\rho$ with the projection $\eT_X\rst_{Y_o}\to\eN_{Y_o/X}$. 
\end{enumerate}
\end{m-lemma}

In general, the mobility conditions are distinct. 
Let $S$, $\YY$ be positive-dimensional, irreducible projective varieties, and consider $\cY=\XX=\YY\times S$. The morphism $\varpi$ is the natural projection, and $\rho$ is the identity. Then~$\cY$ is movable but not generating: any $\cY$-chain has length one. One checks that $\Sigma_2=\emptyset$, reflecting that the members of $\cY$ have no non-trivial intersections.

\begin{proof}
 (\ref{l:mov-1})~  The first three implications are clear. For the last, the surjectivity of $\rd\dta_{o,p}$ at general points $o\in S$, $p\in\YY_o$ implies that any $v\in\eT_{\YY_o,p}\sm\{0\}$ yields a surjective homomorphism $\rd\dta_{o,p}(v)\colon \eT_{S_p,o}\to\eN_{\YY_o/\XX,p}$. Thus, at nearby points $q\in\YY_o$, the evaluation at $q$ is surjective:  
  \[
  H^0(\eI_p\otimes\eN_{\YY_o/\XX})\supset\eT_{S_p,o}\ni\sigma\longmapsto\sigma(q)\in\eN_{\YY_o/\XX,q}.
  \]
  
(\ref{iti:Ygen})~ Each $(\vpi_1,\rho_t)(\cY_t)$ is constructible, and we have 
$S\times\XX={\bigcup}_{t\ges 1} \overline{(\vpi_1,\rho_t)(\cY_t)}$, by hypothesis.  Since $S\times\XX$ cannot be written as a countable union of proper subvarieties, for some $t$, we must have $\overline{(\vpi_1,\rho_t)(\cY_t)}=S\times\XX$. For the converse, there is nothing to prove.

(\ref{l:mov-3})~  The first statement is clear. For strong movability, the infinitesimal deformations of $(\YY_o,p_o)$ within $\cY_o$ are given by pointed subvarieties $(\YY_s,p_s)$, $p_s\in\YY_s\cap\YY_o$, and correspond to $\xi\in\eT_{\cY, p_o}$ such that $\widehat{\rd\rho}_{p_o}(\xi)=0$. The surjectivity of $\widehat{\rd\rho}_{q_o}$ yields  strong movability at the point $((o,o),(p_o,q_o))\in\cY_2$, which projects to $(o,o)\in\Delta_S$. Briefly, $(\rho,\rho)\colon \cY\times_S\cY\to\XX\times\XX$ transversely intersects  $Y_o\times\XX$ at $(p_o,q_o)$. 
\end{proof}

\begin{m-lemma}\label{lm:conn-mov}
Suppose that $\cY$ is such that the self-intersection satisfies $[Y]^2\neq0$. Then we have 
\[
\text{movable}\;\;\Longleftrightarrow\;\;\text{strongly movable.}
\] 
\end{m-lemma}

\begin{proof}
The morphism $\rho\colon \cY\to\XX$ is dominant; by Remark~\ref{rmk:G}(\ref{rmk:G-3}), we have $\cY_2=\cY^2$, $\Sigma_2=S^2$. We deduce that $(\vpi_1,\rho_2)$ is dominant.
\end{proof}

\begin{m-lemma}\label{lm:conn}
Let $\cY$ be a flat family of subvarieties, $D\subset\XX$ an irreducible divisor. The intersection product $[\YY_s]\cdot[D]$ is non-zero, so $\YY_s\cap D\neq\emptyset$ for all $s\in S$ in the following situations:
\begin{enumerate}
\item\label{l:conn-1}
$\XX$ is $\cY$-chain connected in codim-one.
\item\label{iti:Ogen} 
  $\cY$ is generating $\XX$ and $D\cap O_o\neq\emptyset$. Here $o\in S^\gen$ is such that $\YY_o\not\subset D$ and satisfies Lemma~\ref{lm:mov}\,{\rm(\ref{iti:Ygen})}: for some $t\ges 2$, $\rho_t\bigl(\cY|_{\Sigma_t\cap(\{o\}\times S^{t-1})}\bigr)$ contains the open subset $O_o$ of\, $\XX$. 
\item\label{iti:YD} 
$\cY$ is movable and $D\cap\YY_o\neq\emptyset$ for $o\in S^\gen$.
\end{enumerate}

\end{m-lemma}

\begin{proof}
(\ref{l:conn-1})~ Take $o\in S^\gen$ and suppose $[D]\cdot [\YY_o]=0$. The deformation invariance implies that, for all $s\in S$, we have $[D]\cdot[\YY_s]=0$. Then either $D\cap \YY_s=\emptyset$ or $\YY_s\subset D$; otherwise, the intersection $D\cap\YY_s$ is proper, so the intersection product does not vanish.  

 Let $O_o\subset X$ be the open subset swept out by $\cY$-chains starting with $\YY_o$. We have $\codim_\XX(\XX\sm O_o)\ges2$; thus $D\cap O_o\neq\emptyset$.  By chain connectedness, since $D\cap O_o\neq\emptyset$, there are $o=s_0,s_1,\dots,s_n\in S$ such that 
$\YY_{s_{j-1}}\cap \YY_{s_j}\neq\emptyset$ and $D\cap \YY_{s_n}\neq\emptyset$.

 The previous step implies $\YY_{s_n}\subset D$, so $D\cap \YY_{s_{n-1}} \neq\emptyset$, \textit{etc.}; inductively, we deduce $\YY_o\subset D$. Since $o\in S$ is generic, it follows that $\rho(\cY)\subset D$, contradicting that $\rho$ is dominant.

(\ref{iti:Ogen})~ 
 We repeat the argument above for a $\cY$-chain connecting the general points of $\YY_o$ and $D$. If the intersection product vanishes, it follows that $\YY_o\subset D$, and we have a contradiction. 
 
 (\ref{iti:YD})~ 
The intersection is proper, because $\cY$ is movable and $o$ is general in $S$.
\end{proof}

Chain connectedness in codimension one is important; we need conditions ensuring it. 

\begin{m-lemma}\label{lm:>0} 
Let $\cY$ be a flat, movable family of subvarieties of\, $\XX$. 
\begin{enumerate}
\item\label{iti:mori-dream} 
  Let $\XX$ be such that its effective divisors are semi-ample $($e.g.\  $\Pic(\XX)\cong\mbb Z$, or $X$ is a minimal Mori dream space$)$. Suppose $\cY$ is strongly movable. Then any member $\YY$ intersects the effective divisors of\, $\XX$, and $\cY$ is chain connecting in codim-one. 
\item \label{lm:>0-2}
If\, $\cY$ is strongly movable in codim-one, then $\cY$ is chain connecting in codim-one. 
\end{enumerate}
\end{m-lemma}

\begin{proof}
(\ref{iti:mori-dream})~ Suppose there is a $D$ with $[\YY]\cdot[D]\!=\!0$. Since $mD$ is globally generated for some $m\ges1$, $\YY$ is contained in a fibre $F$ of $\XX\srel{\phi}{\to}|mD|$; the same holds for the deformations of $\YY$. Note that $\dim |mD|\ges1$. 
Since $\YY$ is movable, after replacing $\YY$ by a deformation, we may assume that $F$ is a regular fibre of $\phi$. Then $\eN_{F/\XX}\rst_\YY$ (which is trivial, of rank at least one) is a quotient of $\eN_{\YY/\XX}$ (which is $(\dim \YY-1)$-ample, by~\cite{hlc-g2}), and we have a contradiction, so $[\YY]\cdot[D]\!\neq\!0$.

(\ref{lm:>0-2})~ Let $\cY\to\XX$, $S'\subset S$, and $O\subset \XX$ be as in Definition~\ref{def:YY-gener}(\ref{d:YY-g-b})(\ref{d:YY-g-b2}). Take $x, x'\in O$; there are $o,o' \in S'$ such that $x\in\YY_o$, $x'\in\YY_{o'}$. The hypothesis implies that there is a $p\in\YY_o$ such that $\rd\rho_p\colon \eT_{\cY_o,p}\to\eT_{\XX,p}$ is surjective, so $\rho(\cY_o)$ is open about the point $(o,p)$. In particular, the $\cY$-chains of length two starting with $\YY_o$ cover an open subset $U_o\subset\XX$. By repeating the argument, we obtain an open subset $U_{o'}$ swept out by the length two chains starting with $\YY_{o'}$. Since $U_o\cap U_{o'}\neq\emptyset$, there is a $\cY$-chain of length at most four connecting $x$ to~$x'$. 
\end{proof}


\subsection{Irreducible pointed deformation spaces and fundamental groups} \label{ssct:irr} 

\begin{m-definition}\label{def:irr-fibre}
Let $\cY\srel{(\varpi,\rho)}{\longrightarrow}S\times\XX$ be a normal, movable family of irreducible subvarieties. We say that~$\cY$ has {\em irreducible general isotropy} if  $\rho^{-1}(\xx)$ is irreducible for all points $\xx$ in a non-empty open subset $\XX^\iris_\cY\subset\rho(\cY)$. Equivalently, the geometric generic fibre of $\rho$ is connected ($\cY$ is normal). 
For $\xx\in\XX$, $S_\xx:=\rho^{-1}(\xx)=\{ s\mid \xx\in\YY_s \}$ can be viewed as a subset of $S$. The term ``isotropy'' is chosen in analogy with the case of group actions (see Proposition~\ref{prop:irr-fibre}\ref{iti:GOX}). 
\end{m-definition}

\begin{m-lemma}\label{lm:irr-fibre}\leavevmode
\begin{enumerate}
\item \label{l:if-1}
The family $\cY$ has irreducible general isotropy if and only if  $\cY\times_\XX\cY$ possesses exactly one irreducible component dominating $\XX$ $($which contains the diagonal $\Delta_\cY)$.
\item \label{l:if-2}
Suppose $\cY$ has irreducible general isotropy, and let $\XX'\srel{f}{\to}\XX$ be a finite $($surjective$)$ morphism, with $\XX'$ irreducible. Then the following properties hold true: 
\begin{itemize}
\item $\cY\times_\XX\XX'$ contains a unique irreducible component $\cY'$ which dominates $\XX'$. It contains $f^{-1}(\YY_s)$ for any $s\in S$ such that $\YY_s\cap\XX^\iris_\cY\neq\emptyset$.
\item We have $\cY\times_\XX U'\subset\cY'$, where $U'\subset\XX'$ is the \'etale locus of $f$. 
\end{itemize}
\end{enumerate}
\end{m-lemma}

\begin{proof} 
(\ref{l:if-1})~ For the necessity, note that the general fibre of $\cY\times_\XX\cY\to\XX$ is irreducible. For the sufficiency, consider the Stein factorization $\cY\to\hat\cY\to\XX$ (induced by some completion of $\cY$). Should the degree of the second map be greater than one, $\hat\cY\times_\XX\hat\cY$ possesses other components besides the diagonal, determined by pairs of distinct points in the fibres of $\hat\cY\to\XX$. Their pre-image in $\cY\times_\XX\cY$ would still dominate $\XX$.

(\ref{l:if-2})~ The general fibres of $\cY\times_\XX\XX'\to\XX'$ and $\cY\to\XX$ are isomorphic, which implies that 
\[
\cY'=\overline{{\cY'}^{\,\open}}^{\cY\times_\XX\XX'},\quad\text{where } {\cY'}^{\,\open}:=\rho^{-1}(\XX^\iris_\cY)\times_{\XX^\iris_\cY}f^{-1}(\XX^\iris_\cY).
\]
Furthermore, since $\YY_s\cap\XX^\iris_\cY\neq\emptyset$ and $f$ is finite, every component of $f^{-1}(\YY_s)$ dominates $\YY_s$, hence intersects $f^{-1}(\XX^\iris_\cY)$. 

Now proceed to the second claim. (The argument was provided by the anonymous referee.) The \'etale locus $U'\subset\XX'$ of $f$ is Zariski open, irreducible. Then $\cY\times_\XX U'\to \cY$ is \'etale, so it is an open morphism. Thus any component of $\cY\times_\XX U'$ dominates $\cY$ and hence $U'$. Therefore, $\cY\times_\XX U'$ is irreducible and is contained in $\cY'$.  
\end{proof}

Recall, see~\cite[Corollary 5.1]{ver}, that, since $\cY\srel{\rho}{\to}\XX$ is dominant, there is an open $O\subset\XX$  such that $\cY_O:=\cY\times_\XX O\srel{\rho_O}{\to} O$ is topologically locally trivial. One has the following exact sequence in homotopy: 
\begin{equation}\label{esh}
\pi_1^\tpl(\cY_O)\srel{\rho_*}{\longrightarrow}\pi_1^\tpl(O)\lra\pi_0^\tpl(\text{general fibre of}\;\rho)\lra 0.\tag{\rm{esh}}
\end{equation}
Note that, in general, $\XX\sm O$ contains divisors. The general fibre is connected if $\rho_*$ is surjective. 

\begin{m-proposition}\label{prop:irr-fibre} 
The general isotropy of\, $\cY$ is irreducible in any of the situations below:  
\begin{enumerate}
\item\label{iti:GOX} 
A group $G$ acts on $\XX$ with an open orbit $O$ and $\xx_0\in O$ has connected stabilizer $($isotropy group$)$. Here we let $\cY=G\times\YY$, where $\YY\subset\XX$ is irreducible and intersects $O$. 

\item\label{iti:cod2} 
\begin{enumerate}
\item\label{cod2-a}
  The image $\rho(\cY)$ contains an open subset $O$ with  $\codim_\XX(\XX\sm O)\ges2$, and the smooth locus of $\rho$ is dense in $\rho^{-1}(x)$, for all $x\in O$ $($e.g.~$\cY_O\srel{\rho_O}{\to} O$ is smooth$)$;
\item \label{cod2-b}
there exist a normal, projective $\widebar\cY\supset\cY$ and a morphism $\widebar\cY\srel{\widebar\rho}{\to} X$ extending $\rho$ such that  
\[
\forall\,x\in O,\;\;\rho^{-1}(x)\subset{\widebar\rho}^{-1}(x) \text{ induces a bijection of the connected components;} 
\]
\item \label{cod2-c}
$\XX$ is algebraically simply connected $($e.g.~$\pi_1^\tpl(\XX)=\{1\})$.
\end{enumerate}
\item \label{iti:esh}  
The image $\rho(\cY)$ contains an open subset $O$ with $\codim_\XX(\XX\sm O)\ges2$, with the following properties: 
\begin{itemize}
\item[$\varstar$]
There is a smooth member $\YY\in\cY$ such that $\codim_\YY(\YY\sm O)\ges2$, and 
\begin{equation}\label{pi}
\pi_1^\tpl(\YY)\lra\pi_1^\tpl(\XX) \text{ is surjective.}\tag{$\pi$}
\end{equation} 
\item[$\diamond$] \label{diamond}
$\cY_O\to O$ is locally trivial for the ${\euf C}^0$/Euclidean topology. 
\end{itemize}
\end{enumerate}
\end{m-proposition}

\begin{proof}
(\ref{iti:GOX})~  Let $H_0$ be the stabilizer of $\xx_0$. Then $\rho^{-1}(\xx)=\{ (g,y)\in G\times(\YY\cap O) \mid g\cdot y = \xx \}$ is a smooth fibre bundle over $\YY\cap O$, with fibre isomorphic to $H_0$. 

(\ref{iti:cod2})~  Let $\widebar\cY\srel{g}{\to}\hat\cY\srel{f}{\to}X$ be the Stein factorization of $\widebar\rho$. Suppose the ramification locus of $f$, $R\subset\hat\cY$, is non-empty; then $R$ is a divisor (purity of the branch locus). Assumption (\ref{cod2-b}) implies $f^{-1}(O)=g(\cY_O)$, so~$R$ intersects $g(\cY_O)$, even its smooth locus (normality of $\hat\cY$). Let $\hat y\in R\cap g(\cY_O)$ and $x=f(\hat y)$. Since $f^{-1}(x)$ has multiplicity at $\hat y$, the fibre $\rho^{-1}(x)$ contains a non-reduced component, contradicting (\ref{cod2-a}). Thus $f$ is unramified, so $\hat\cY\to\XX$ is \'etale. We conclude that $\hat\cY=\XX$, by (\ref{cod2-c}), so the general fibre of $\rho$ is irreducible. 

(\ref{iti:esh})~  By the codimension two hypothesis,  $\pi_1^\tpl(O)\to\pi_1^\tpl(\XX)$  and $\pi_1^\tpl(\YY\cap O)\to\pi_1^\tpl(\YY)$  are isomorphisms. Property $(\pi)$ yields the surjectivity of the composed homomorphism $\pi_1^\tpl(\YY\cap O)\to\pi_1^\tpl(\cY_O)\to\pi_1^\tpl(O)$, so the second arrow is surjective too. Finally, the exact sequence (\ref{esh}) yields the connectedness of the general fibre of~$\rho$.
\end{proof}

Condition ``\ref{diamond}''  is restrictive; $\XX\sm O$ should contain no divisors. But the proof uses only the exact sequence (\ref{esh}), for which it is enough that $\rho_O$ is a $2$-quasi-fibration/2-equivalence (see~\cite{may} and references therein). This condition is much weaker; it is satisfied in numerous situations. 

The author~\cite[Corollary 2.3]{hlc-pas} proved Property (\ref{pi}) for the \'etale fundamental groups as soon as $\cd(\XX\sm\YY)\les\dim\XX-2$. Here we point out that this matter, combined with the results~\cite{nap-ram}, yields Property~(\ref{pi}) for partially ample subvarieties. 

\begin{m-theorem}\label{thm:pi}
Let $\YY\srel{\iota}{\to}\XX$ be a lci subvariety of\, $\XX$. We have the following implications: 
\begin{enumerate}
\item\label{thm:pi-1} If\, $\YY$ is strongly movable, then $\iota_*\pi_1^\tpl(\YY)$ has finite index in $\pi_1^\tpl(\XX)$.
\item\label{thm:pi-2} If\, $\YY$ is $1^\pos$ in $\XX$ \emph{(no mobility required)}, then $\pi_1^\tpl(\YY)\srel{\iota_*}{\to}\pi_1^\tpl(\XX)$ is surjective.
\end{enumerate}
\end{m-theorem}

\begin{proof}
In both cases, the normal bundle of $\YY$ is locally free and $(\dim\YY-1)$-ample (\textit{cf.}~\cite[Proposition~3.5]{hlc-g2}, resp.~\cite[Proposition 1.6]{hlc-pas}). Thus, for any line bundle $\eL$ on $\XX$, the vector space $H^0(\hat\XX_\YY,\hat\eL)$ is finite dimensional (\textit{cf.}~\cite[Theorem 2.1(i)]{hlc-g2}, also the remark after~\cite[Theorem 0.1]{nap-ram}). Then~\cite[Theorem~1.1]{nap-ram} implies that  $Q:=\pi_1^\tpl(\XX)/\iota_*\pi_1^\tpl(\YY)$ is finite. 

We continue with (\ref{thm:pi-2}). Let $\tld\XX\srel{f}{\to}\XX$ be the finite cover corresponding to $Q$. Then $\iota$ can be lifted to $\tld\XX$, and the components of $f^{-1}(\YY)$ bijectively correspond to the elements of $Q$. However, the $1^\pos$-property implies $\cd(\XX\sm\YY)\les\dim\XX-2$, so $f^{-1}(\YY)$ is connected (\textit{cf.}~\cite[Theorem 2.2]{hlc-pas}). It follows that $Q=\{1\}$, so $\iota_*$ is surjective.
\end{proof}

Thus we have the implication $1^\pos\Rightarrow(\pi)$. Corollary~\ref{cor:strg-mov} shows that the converse $(\pi)\Rightarrow 1^\pos$ holds under mobility assumptions on $\YY$.


\section{A Chow-type G3-criterion and applications}\label{sct:chow}

\subsection{The criterion} 
We prove that mobility conditions imply the G3-property. This result is inspired from~\cite{chow} and~\cite[Theorem 13.4]{bad}, which consider subvarieties of homogeneous varieties.

\begin{m-theorem} \label{thm:g3}
Let $\YY_o$ be a G2-subvariety of\, $\XX$ which has the following properties: 
\begin{itemize}[leftmargin=3ex]
\item It belongs to a flat, movable family $\cY={(\YY_s)}_{s\in S}$ with irreducible general isotropy. 
\item Any effective divisor of\, $\XX$ intersects $Y_s$, $s\in S^\gen$. $($Equivalently, by Lemma~\ref{lm:conn}\,{\rm(\ref{iti:YD})}, the intersection product of $\YY_o$ with any effective divisor of\, $\XX$ is numerically non-trivial.$)$
\end{itemize}
Then $\YY_o$ is G3 and $\cd(\XX\sm\YY_o)\les\dim\XX-2$ in either one of the following situations:
\begin{enumerate}
\item\label{thm:g3-1} $\XX$ is algebraically simply connected $(\pi_1^\etl(\XX)=\{1\})$.
\item\label{thm:g3-2} $[\YY_o]^2\neq0\in A^*(\XX)/\sim_{\num}$.  
\end{enumerate}
\end{m-theorem}

\begin{proof}
Since $Y_o$ is G2, there exist a normal, irreducible, projective variety $\XX'$, a G3-subvariety $\YY'_o\subset \XX'$, and a finite, surjective morphism $f\colon (\XX',\YY'_o)\to(\XX,\YY_o)$ such that $f$ is \'etale in a neighbourhood of $\YY'_o$ (\textit{cf.}~\cite[Corollary 9.20]{bad}). 
Let $D'\subset \XX'$ be the ramification divisor of $f$ and $D:=f_*(D')$; we want $D=\emptyset$, that is, $f$ is \'etale. 
Let us consider the diagram  
\[
\xymatrix@R=1.89em@C=2em{
\cY'_{\ovl{\omega}}\ar[r]\ar[d]^-{\vphi_{\ovl{\omega}}}&
\cY'_\omega\ar[r]\ar[d]^-{\vphi_\omega}
&\cY'\subset\cY\times_\XX\XX'\ar[r]^-{\rho'}\ar[d]^-\vphi&\XX'\ar[d]^-f
\\ 
\cY_{\ovl{\omega}}\ar[r]\ar[d]&\cY_\omega\ar[r]\ar[d]&\cY\ar[r]^-\rho\ar[d]^-\vpi&\XX\rlap{.}
\\ 
\ovl{\omega}:=\Spec(\ovl{\kk(S)})\ar[r]&{\omega}:=\Spec({\kk(S)})\ar[r]&S&
}
\]
Here $\cY'\subset\cY\times_\XX\cY$ is the irreducible component introduced in Lemma~\ref{lm:irr-fibre}, dominating $\XX'$. 
The ramification divisor of $\vphi\colon \cY'\to\cY$ is $D_1:=(\rho')^{-1}(D')\cap\cY'$. We are going to show that $D_1\cap \cY'_\omega=\emptyset$; it suffices to prove that $D_1\cap \cY'_{\ovl{\omega}}=\emptyset$, that is, $\vphi_{\ovl{\omega}}$ is \'etale. 

Since $f$ is \'etale about $\YY'_o$, Lemma~\ref{lm:irr-fibre} implies that $\YY'_o=\YY_o\times_\XX \YY'_o\subset\cY'$. Thus $\YY'_o$ is a component of $\cY'_o=(\vpi\vphi)^{-1}(o)$ and $\vphi$ is \'etale in its neighbourhood. But $\YY'_o$ is the specialization of a component $Z$ of the fibre $\cY'_{\ovl{\omega}}\to{\ovl{\omega}}$, so $\vphi_{\ovl{\omega}}\colon Z\to\cY_{\ovl{\omega}}$ is \'etale (because $\vphi$ is so about $\YY'_o$). The components of $\cY'_{\ovl{\omega}}$ are conjugate under the Galois group ${\rm Gal}\big(\ovl{\kk(S)}/\kk(S)\big)$, since $\cY'_\omega$ is irreducible. Thus $\vphi_{\ovl{\omega}}$ is \'etale everywhere and so is~$\vphi_\omega$. This proves that $D_1\cap\cY'_\omega=\emptyset$. 

Lemma~\ref{lm:irr-fibre} says that  $\cY'$ contains $f^{-1}(\YY_s)$ for $s\in S^\gen$. We deduce that $D_1\cap f^{-1}(\YY_s)=\emptyset$. Hence $D\cap\YY_s=\emptyset$; the subvariety $\YY_s$ does not intersect the branch divisor. This contradicts that $\YY_s\cap D\neq\emptyset$ for $s\in S^\gen$ (or that $[\YY_o]\cdot[D]\neq0$ ) unless $D=\emptyset$, so $f$ is \'etale.

(\ref{thm:g3-1})~ If $\XX$ is algebraically simply connected, then $f$ is an isomorphism, so $\YY_o$ is G3 in $\XX$. 

(\ref{thm:g3-2})~ Now suppose $d:=\deg(f)\ges2$, and let $\YY_o',\YY_1',\dots,\YY_{d-1}'$ be the pre-images of $\YY_o$ by $f$; they are deformation equivalent.  Then $[\YY_o]^2$ can be computed on the covering space $\XX'$: 
$$
d\cdot[\YY_o]^2=[\YY_o']\cdot [f^{-1}(\YY_o)]
=[\YY_o']\cdot([\YY_o']+[\YY_1']+\dots+[\YY_{d-1}'])
=d\cdot[\YY_o']\cdot[\YY_1']=0.
$$
\nit This contradicts the hypothesis, so $d=1$ and $f$ is \'etale.
In both situations, $Y_o$ is G3 and intersects all effective divisors, so Speiser's criterion implies that $\cd(\XX\sm\YY_o)\les\dim\XX-2$. 
\end{proof}

\nit Note that, if $[\YY_o]^2\neq0$, then it is strongly movable (\textit{cf.}~Lemma~\ref{lm:conn-mov}). For lci subvarieties, strong mobility yields the G2-property (see Section~\ref{ssct:setup}). 

\begin{m-corollary}\label{cor:mori}
Let $\YY\subset\XX$ be lci, strongly movable, belong to a flat family $\cY$ with irreducible general isotropy, and intersect numerically non-trivially every effective divisor of\, $\XX$. Then $\YY$ is G3, actually $1^\pos$, in either one of the following cases: 

\begin{enumerate*}
\item\label{case1}  $\XX$ is rationally connected; \;\unbar{\emph{or}}\item\label{case2} $[\YY]^2\neq0$  $($\textit{e.g.}~$\eN_{\YY/\XX}$ is ample, $2\dim\YY\ges\dim\XX)$.\end{enumerate*}  
\end{m-corollary}

\begin{m-corollary}\label{cor:dbar}
Suppose $\YY$ is lci, strongly movable, belongs to a flat family $\cY$ with irreducible general isotropy, and intersects numerically non-trivially all effective divisors of\, $\XX$. Then there exist an \'etale cover $\XX'$ of\, $\XX$ and an inclusion $\YY'\subset\XX$, such that $\YY'\to\YY$ is an isomorphism and $\YY'$ is G3 $($actually $1^\pos)$ in $\XX'$.
\end{m-corollary}

This answers, in a broader context, a question raised by Debarre~\cite[Conjecture 2.9]{dbar-mumford} about generating subvarieties of homogeneous varieties, which are known to intersect all divisors. The isotropy condition is satisfied in this situation, too, by Proposition~\ref{prop:irr-fibre}(\ref{iti:GOX}).

\begin{proof}
The proof of Theorem~\ref{thm:g3} shows that $f$ is the desired \'etale morphism. (We stop before discussing cases (\ref{case1}) and (\ref{case2}).)
\end{proof}

\begin{m-remark}\label{rk:gold}
Suppose $\YY\subset\XX$ is smooth. Then $\eN_{\YY/\XX}$ is a quotient of $\eT_\XX$, so the co-normal bundle is automatically not pseudo-effective as soon as $\eT_\XX$  is $(\dim\YY-1)$-ample. 

Let $\XX=G/P$ be rational homogeneous, with $G$ semi-simple. Goldstein~\cite{gold} computed the partial amplitude $\amp(\XX)$ of the tangent bundle; it satisfies $\amp(X)\les\dim\XX-\ell$, where $\ell$ is the minimal rank of the simple factors of $G$. 
Thus the normal bundle of $\YY$ is $(\dim Y-1)$-ample --~not pseudo-effective~-- as soon as $\dim\XX-\ell\les\dim\YY-1$, that is, $\codim_\XX\YY<\ell$. In this situation, $\YY$ also intersects non-trivially every effective divisor of $\XX$ (\textit{cf.} \cite[Lemma 13.17]{bad}). 

Theorem~\ref{thm:g3} implies that $\YY$ is G3 in $\XX$; thus we recover Faltings' result~\cite[Satz 8]{falt-homog}; see also~\cite[Corollary 13.18]{bad}. Although here we require $\YY$ to be smooth, our approach places the G3-property in a much broader context.
\end{m-remark}


\subsection{Examples}

The previous results show the importance of mobility for the $G3$-property. Let $\cY\srel{(\vpi,\rho)}{\subset} S\times \XX$ be a family of lci subvarieties of $\XX$. We are going to discuss two fairly independent situations where Theorem~\ref{thm:g3} applies:
\begin{itemize}
\item Every effective divisor of $\XX$ is semi-ample (some multiple is globally generated).
\item The family $\cY$ is strongly movable.
\end{itemize}


\subsubsection{Minimal Mori dream spaces}\label{sssct:mds}
They are the prototype of varieties whose effective divisors are semi-ample (\textit{cf.}~\cite[Proposition 1.11]{hu-keel}). Examples include numerous Fano, toric, and spherical varieties, and GIT-quotients. 

\begin{m-corollary}\label{cor:mori-g3} 
Suppose that every effective divisor of\, $\XX$ is semi-ample. Then the lci subvariety $\YY$ is $1^\pos$ in $\XX$,  thus G3, in the following situation:
\begin{enumerate}
\item $\YY$ is strongly movable $($\textit{e.g.}~very moving$)$,  
and $\cY$ has irreducible general isotropy; 
\item \unbar{\emph{either}} $\pi_1^\etl(\XX)=\{1\}$, \unbar{\emph{or}} $\YY$ has ample normal bundle and $2\dim\YY\ges\dim\XX$.
\end{enumerate}
\end{m-corollary}

\begin{proof} 
Lemma~\ref{lm:>0} implies that $\YY$ intersects every effective divisor of $\XX$ non-trivially. Finally, in the case $\pi_1^\etl(\XX)\neq\{1\}$, the excess intersection theorem implies $[Y]^2\neq0$. 
\end{proof}


\subsubsection{Strongly movable families}\label{sssct:strg-mov}

\begin{m-corollary}\label{cor:strg-mov}
Let $\cY$ be a family of lci subvarieties of\, $\XX$,  strongly movable  in codim-one. Let $\YY$ be a strongly movable member of\, $\cY$. We assume moreover 
\vskip.5ex
\begin{tabular}{rl}
\unbar{\emph{either}} {\rm(i)}\,& $\pi_1^\etl(\XX)=\{1\}$, \textit{e.g.}~$X$ is rationally connected,  and $(\cY,\rho)$ satisfies  Conditions~\ref{prop:irr-fibre}\,{\rm{(\ref{iti:cod2})}};
\\ 
\unbar{\emph{or}} {\rm(ii)}& $[\YY]^2\neq0$ and $(\cY,\rho)$ satisfies Conditions~\ref{prop:irr-fibre}\,{\rm(\ref{iti:esh})}. 
\end{tabular}

\noindent Then $\YY$ is G3 $($actually $1^\pos)$ in $\XX$.
\end{m-corollary}

\begin{proof}  
The strong movability implies that $\YY$ is G2 in $\XX$, and it intersects all the effective divisors, by Lemma~\ref{lm:>0}. The general isotropy is irreducible by Proposition~\ref{prop:irr-fibre}. 
\end{proof}


\subsubsection{Subvarieties of (almost) homogeneous varieties}\label{ssct:subvar-homog}

A common situation ensuring the mobility of subvarieties is when the ambient variety possesses a large group of automorphisms. Let $G$ be a connected linear algebraic group with identity $e$, and consider the morphism  
$$
\gamma\colon G\times G\lra G,\quad\gamma(g',g):=g'g^{-1}.
$$
Let $\XX$ be a smooth almost homogeneous $G$-variety with open orbit $O$. Note that $\XX$ is automatically rationally connected. 
We set $\partial\XX:=\XX\sm O$; it is the finite union of its irreducible components. Let $\{\partial\XX\}_{\rmdiv}$ be the set of one-codimensional components, called \emph{boundary divisors}. All the components of $\partial\XX$ are $G$-invariant; if $G$ acts with finitely many orbits, they are closures of $G$-orbits. Examples are varieties with a finite number of orbits for the action of a linear reductive group. 

Let $H$ be the stabilizer of a point $y_o\in O$, so $O\cong G/H$; we assume that $H$ is connected. Let $G\times\XX\srel{\mu}{\to} \XX$ be the action and, for $x\in \XX$, let $\mu_x(\cdot):=\mu(\cdot,x)$. 
Suppose $\YY\subset \XX$ is an irreducible subvariety intersecting $O$; we set 
$\YY_O:=\YY\cap O$ and $\partial\YY:=\YY\sm\YY_O$. 
We retrieve the situation~\eqref{eq:syx}:
\begin{equation}
\xymatrix@R=0.99em{
&\cY:=G\times \YY\ar[dr]_-\mu\ar[dl]^-\vpi&
\\  
G&&\XX 
}\label{eq:strg-gen}
\end{equation}
This family has irreducible general isotropy, by Proposition~\ref{prop:irr-fibre}(\ref{iti:GOX}). 

\begin{m-definition}\label{def:gener} 
We consider the following objects: 
\begin{itemize}
\item For $y_o\in \YY_O$, let $G_{\YY,y_o}:=\mu^{-1}_{y_o}(\YY)$; it is a closed subvariety of $G$. Denote by $G_\YY$ the subgroup generated by $G_{\YY,y_o}$; it is closed in $G$, independent of $y_o\in \YY_O$.

\item Define $S_\YY:=\gamma(G_{\YY,y_o},G_{\YY,y_o})$; it is constructible in $G$, contains $G_{\YY,y_o}$, and is independent of $y_o\in \YY_O$. In fact, $S_\YY$ consists of those elements of $G$ which send some point of $\YY_O$ to another point of $\YY_O$; we recover $\Sigma_o$ in Diagram~\eqref{eq:cY}, for $o=e\in G$.
\end{itemize}
We say that $\YY$ \emph{generates} $\XX$ if $G_{\YY}=G$; it \emph{strongly generates} $\XX$ if $\mu(S_\YY,\YY_O)$ is open. These, respectively, correspond to the generating and strongly movable properties in Definition~\ref{def:YY-gener}.
\end{m-definition}

Note that, for homogeneous $\XX=G/P$, the subset $S_\YY$ is closed in $G$. Indeed, if $y_o=e$, then $P\subset G_{\YY,y_o}$ and~$\gamma$ factorizes through $(G\times G)/P$, where $P$ acts diagonally. But the morphism $(G\times G)/P\to G$ is projective, so the image of $(G_{\YY,y_o}\times G_{\YY,y_o})/P$ is closed in $G$. Let us remark that the strong generation property is satisfied if $G$ acts generically $2$-transitively on $\XX$ (that is, the diagonal action has an open orbit in $\tld O\subset\XX\times\XX$, see~\cite{pop}) and $\YY\times\YY$ intersects $\tld O$. 

\begin{m-lemma}\label{lm:strg-gen}
Suppose $\YY$ strongly generates $\XX$. Then the following properties hold:
\begin{enumerate}
\item\label{l:sg-1} The subvariety $\YY$ generates $\XX$.
\item\label{l:sg-2} The family $\cY=G\times\YY$ in Diagram~\eqref{eq:strg-gen} is strongly movable.
\item\label{l:sg-3} If $\codim_\XX(\partial\XX)\ges2$, then $\XX$ is $\cY$-chain connected in codim-one.
\item\label{iti:DO} For any irreducible divisor $D\subset\XX$ with $D\cap O\neq\emptyset$, we have $[\YY]\cdot [D]\neq0$.
\end{enumerate}
\end{m-lemma}

\begin{proof}
(\ref{l:sg-1})~ The stabilizer $H\subset G$ of $y_o$ is contained in $G_{\YY,y_o}$, and 
$$
\mu(S_\YY,\YY_O)=\mu\big((G_{\YY,y_0}\cdot G_{\YY,y_o}^{-1}\cdot G_{\YY,y_o}),y_o\big)
\cong 
(G_{\YY,y_0}\cdot G_{\YY,y_o}^{-1}\cdot G_{\YY,y_o})/H.
$$
\nit If $\mu(S_\YY,\YY_O)$ is open in $O\cong G/H$, then $G_{\YY,y_0}\cdot G_{\YY,y_o}^{-1}\cdot G_{\YY,y_o}$ is $H$-invariant and open in $G$. Thus $G_\YY\subset G$ is closed and contains an open subset; hence it coincides with $G$. 

(\ref{l:sg-2})~ For $g\in G$, we have $S_{g\YY}=gS_\YY g^{-1}$ and $\mu(S_{g\YY},g\YY_O)=g\mu(S_\YY,\YY_O)$. The right-hand side is open in $\XX$; hence the morphism $(\vpi_1,\rho_2)$ in Definition~\ref{def:YY-gener} is dominant.

(\ref{l:sg-3})~ Part~(\ref{l:sg-1}) implies that $S_{\YY}$ generates $G$: for any $g\in G$, there are $e=g_0,g_1,\dots,g_n\in S_{\YY}$ whose product is~$g$. Then 
$\YY_0:=\YY, \YY_1:=g_1\YY,\dots,\YY_n:=g_1\dots g_n\YY$ connects $y_o$ to $gy_o$. 

(\ref{iti:DO})~ The family $\cY$ generates $\XX$, so can we apply Lemma~\ref{lm:conn}(\ref{iti:Ogen}) with $t=2$.
\end{proof}

\begin{m-theorem}\label{thm:G3-a-homog} 
Let $\YY\subset \XX$ be an lci subvariety with $\YY_O\neq\emptyset$,  with the following properties: 
\begin{enumerate}
\item\label{t:G3ah-1} It strongly generates $\XX$.
\item\label{t:G3ah-2} Its intersection with every boundary divisor is non-empty $($e.g.~$\codim_\XX(\partial\XX)\ges 2)$.
\end{enumerate} 
Then $\YY$ is $1^\pos$; in particular, it is G3 in $\XX$.  
\end{m-theorem}

\begin{proof}
Since $G$ is linear, $\XX$ is simply connected. By Lemma~\ref{lm:strg-gen}, $\YY$ is G2. To apply Theorem~\ref{thm:g3} to $\cY=G\times\YY$, we should verify that any effective divisor $D$ of $\XX$ intersects the general translates of $\YY$. If $D\in\{\partial\XX\}_{\rmdiv}$, then it is $G$-invariant and $\YY$ intersects it, by hypothesis. If $D\cap O\neq\emptyset$, Lemma~\ref{lm:strg-gen}(\ref{iti:DO}) says that $[\YY]\cdot [D]\neq0$.
\end{proof}

For homogeneous varieties, the G3-property of the diagonal in $\XX \times \XX$ is proved in~\cite[Theorem~13.19]{bad}). The corollary below shows that our techniques yield this property in many, previously unknown, situations (\textit{e.g.}~$X$ is a compactification of $G/T$, where $T$ is a maximal torus in the reductive group $G$).  

\begin{m-corollary}\label{cor:DXX}
Suppose the stabilizer of some $($any$)$ point $y_o\in O$ contains a Cartan subgroup of\, $G$. 
Then the diagonal is $1^\pos$, thus G3, in the product $\XX\times \XX$. 
\end{m-corollary}

\begin{proof}
We apply Theorem~\ref{thm:G3-a-homog} to the diagonal $\Delta_\XX\subset \XX\times \XX$. The group $G\times G$ acts on $\XX\times \XX$ with open orbit $O\times O$. We verify that $\Delta_\XX$ is strongly generating, that is, $S_{\Delta_\XX}\cdot O$ is open in $O\times O$. Let $H$ be the stabilizer of the point $y_o\in O$. A direct computation yields  
\[
S_{\Delta_\XX}
=\{ (g,g)\cdot(e,\Ad_a(h))\mid a,g\in G,h\in H\}\quad(\text{$e\in G$ is the identity}). 
\]
\nit Since $H$ contains a Cartan subgroup of $G$, ${\bigcup}_{a\in G}\Ad_a(H)$ contains an open subset of $G$ (\textit{cf.} \cite[Theorem~11.10]{bor}), so $S_{\Delta_\XX}$ contains an open subset of $G\times G$. 

Let $\tld D\subset\XX\times\XX$ be an irreducible divisor; we verify that it intersects the general translates of $\Delta_\XX$. If it is contained in the boundary of $\XX\times\XX$, then $\tld D=D\times\XX$ (or $\XX\times D$), with $D\in\{\partial\XX\}_{\rmdiv}$, and $\tld D$ intersects $\Delta_\XX$. Otherwise, Lemma~\ref{lm:strg-gen}(\ref{iti:DO}) implies $[\Delta_\XX]\cdot[\tld D]\neq0$. 
\end{proof}


Let $\widebar G$ be an equivariant completion of $G$. Let $\VV\subset\XX$ be a projective subvariety which intersects the open orbit $O$, so $\cal\VV:={G\times\VV}\srel{\mu}{\to}\XX$ is dominant. 
We identify $\cal\VV$ with the graph of $\mu$, so $\cal\VV=\{ (g,v,x) \mid \mu(g,v)=x \}$, and take its closure $\widebar{\cal V}$ in $\widebar G\times\VV\times\XX$. It admits morphisms 
$$
\xymatrix@R=0.99em@C=3.5em{ &\widebar{\cal\VV}\ar[dr]\ar[dl]^-{\vpi}\ar[drr]^-{\widebar\mu}&\\ 
\widebar G&&\VV\times\XX\; \ar[r]&X.}
$$

\begin{m-corollary}\label{cor:barV}
Let $\YY\subset\XX$ be as in Theorem~\ref{thm:G3-a-homog}. Then ${\widebar\mu}^{-1}(\YY)$ is connected, actually G3. 
\end{m-corollary}

The statement is similar in flavour to~\cite[Th\'eor\`eme 2.4(a)]{dbar-mumford}.

\begin{proof}
Note that $\widebar{\mu}\colon \widebar{\cal\VV}\to\XX$ is surjective. Since $\YY$ is G3 in $\XX$, Hironaka--Matsumura's result~\cite{hir-mat} implies that the pre-image is still G3; in particular, it is connected.
\end{proof}

The connectedness of $\VV\cap\YY$  itself requires stronger assumptions. For a $G$-orbit $O'$ of $\XX$, let $\VV_{O'}:=\VV\cap O'$; we use similar notation for $\YY$. For $(v,x)\in\VV\times\XX$, both in the same $G$-orbit of $\XX$, we define $G_{v,x}:=\{ g\in G\mid \mu(g,v)=x \}$. 
The following is an analogue of Mumford's theorem, see~\cite[Th\'eor\`eme 1.1]{dbar-mumford}, for homogeneous varieties. The extra hypotheses~(\ref{iti:div}) ensure that $\VV$ intersects effective divisors of $\YY$ and \textit{vice versa}; this condition is assumed in \lcit 

\begin{m-theorem}\label{thm:mum}
Let $G$ be simply connected, acting with finitely many orbits on $\XX$. Let $\VV,\YY$ be irreducible subvarieties with $\dim\VV+\dim\YY>\dim\XX$, satisfying the properties below: 
\begin{enumerate}
\item\label{iti:VY} The $G$-orbits intersect $\VV$ and $\YY\!$, and also their singular loci, in the expected dimension. 

\nit More precisely: for each stratum $O'\subset\XX$,  the locus 
$$
\VV_{O'}^{\rm smooth}:=
\bigl\{v\in\VV_{O'}\mid v\in\VV^{\rm smooth}\cap O',\;\text{the intersection is transverse at}\;v\bigr\}
$$
is open and dense in $\VV_{O'}$; a similar condition holds for $\YY_{O'}$.
\item\label{iti:p01} For each stratum $O'\subset\XX$ and $(v,y)\in(\VV_{O'}\times\YY_{O'})^\gen$, the following holds:
\begin{m-eqn}{
    \codim_{G_{v,y}}\left\{ g\in G_{v,y}\mid g_*(\eT_{\VV,v})+\eT_{\YY,y}\neq\eT_{\XX,y} \right\}>
 \begin{cases} 
   1&\text{if }\codim_\XX(O')=0,1,\\[.2ex]
   0&\text{if }\codim_\XX(O')>1.
\end{cases}
}\label{eq:p1}
\end{m-eqn} 
$($These are Mumford's $P(1)_{O',(v,y)}$- and $P(0)_{O',(v,y)}$-conditions.$)$ 
\item\label{iti:div} Either one of the following two assumptions holds: 
\begin{enumerate}
\item\label{itii:amp} Both $\VV$ and $\YY$ are lci and have ample normal bundle. 
\item\label{itii:int} The $P(0)_{O,(v,v',y)}$-condition holds for all $v, v'\in\VV_{O}^\gen$, $y\in\YY_{O}$; similarly, the $P(0)_{O,(y,y',v)}$-condition holds for all $y, y'\in\YY_{O}^\gen$, $v\in\VV_{O}$:
\begin{m-eqn}{
\begin{scalebox}{.975}
{$\codim_{G_{v,y}\times G_{v',y}}\left\{ (g,g')\!\mid\! (g,g')_*\eT_{\VV\times\VV,(v,v')}+\Delta_*\eT_{\YY,y}\neq\eT_{\XX\times\XX,(y,y)} \right\}>0.$}
\end{scalebox}
}\label{eq:p0}
\end{m-eqn} 
Here $\YY\srel{\Delta}{\to}\XX\times\XX$ is the diagonal map. 
\end{enumerate}
\end{enumerate}
Then $\VV\!\cap\YY$ is non-empty, connected. For $g\in G^\gen$, $g\VV\!\cap\YY$ is irreducible of the expected dimension. 
\end{m-theorem}

\begin{proof}
\emph{Claim}~1.~ $\VV$ intersects non-trivially every divisor of $\YY$ and \textit{vice versa}. 

Let $D_\YY\subset\YY$ be an irreducible divisor. If $D_\YY\subset\YY\cap\partial\XX$, then~(\ref{iti:VY}) implies that $D_\YY$ is a component of~$\overline{\YY_\BD}$, where $\BD$ is a one-codimensional orbit ($G$ has finitely many orbits, so boundary divisors are closures of orbits). But we know $\VV_\BD\neq\emptyset$, by the same~(\ref{iti:VY}); thus a translate of $\VV$ will intersect $D_\YY$. In the case $D_\YY\cap\YY_O\neq\emptyset$, we reach the same conclusion. 

Let us consider situation~(\ref{itii:amp}). Since a translate $g\VV$ intersects $D_\YY$, the excess intersection theorem implies that  $[\VV]\cdot[D_\YY]=[g\VV]\cdot[D_\YY]\neq0$. 

Next we analyse situation~(\ref{itii:int}). 

\emph{Suppose $D_\YY\cap\YY_O\neq\emptyset$.}~ Condition~\eqref{eq:p0} with an equality sign means that $g\VV\times g'\VV$ transversally intersects $\Delta_\YY\subset\XX^2$; equivalently, $g\VV\cap\YY$ and $g'\VV\cap\YY$ are transverse in $\YY$. Let $p\in D_\YY\cap\YY_O$ be a smooth point. There are $g_p, g'\in G$, $v_p, v'\in\VV_O$ as in~(\ref{iti:p01}) and~(\ref{iti:div})(\ref{itii:int}), such that the intersections $g_{p}\VV\cap\YY$, $g'\VV\cap\YY$, $g_{p}\VV\cap g'\VV\cap\YY$ are transverse at $y$: 
\begin{align*}
&g_{p*}\eT_{\VV,v_p}+\eT_{\YY,p}=\eT_{\XX,p},\;\;g'_*\eT_{\VV,v'}+\eT_{\YY,p}=\eT_{\XX,p},\\
&(g_{p*}\eT_{\VV,v_p}\cap\eT_{\YY,p})+(g'_*\eT_{\VV,v'}\cap\eT_{\YY,p})=\eT_{\YY,p}.
\end{align*}
Then at least one of $g_p\VV\cap\YY$ or $g'\VV\cap\YY$ intersects $D_\YY$ transversally at $p$: both tangent spaces cannot be contained in $\eT_{D_\YY,p}$. (What follows is basically the proof of~\cite[Theor\`eme 1.1(1)]{dbar-mumford}.) 

Assume the former case.
Then $G\times\VV\srel{\mu}{\to}\XX$ is transverse to $D_\YY$ at $(g_p, v_p)$. Since $v_p$ belongs to $O$, $\mu$ is open at $(g_p, v_p)$, so $\mu^{-1}(D_\YY)\subset G\times\VV$ is smooth at $(g_p, v_p)$; its local dimension is $\dim H+\dim\VV+\dim\YY-1\ges\dim G$. The transversality implies that the differential of $\mu^{-1}(D_\YY)\xrightarrow{\vpi_D}G$ is surjective at $(g_p, v_p)$; thus $\vpi_D(\mu^{-1}(D_\YY))$ contains an open subset of $G$. But $\vpi_D$ is proper, so $\mu^{-1}(D_\YY)\xrightarrow{\vpi_D}G$ is surjective, which implies $\vpi_D^{-1}(e)=\VV\cap D_\YY\neq\emptyset$.\hfill$\Box$

\emph{Suppose $D_\YY\cap\YY_O=\emptyset$.}~ 
So $D_\YY$ is a component of $\overline{\YY_\BD}$, where $\BD$ is a one-codimensional orbit. We repeat the reasoning above for $\mu\colon G\times\overline{V_\BD}\to \overline\BD$. (Note: we use condition~(\ref{iti:p01}); if $g_p\VV$, $\YY$ are transverse at $p\in\BD$, then so are  $g_p\VV_\BD, \YY_\BD$.)  Hence the components of $\YY_\BD$ (that is, $D_\YY$) intersects every $G$-translate of the components of $\VV_\BD$. \hfill$\Box$\smallskip 

Let $Z\!:=(G\times\VV)\times_\XX\YY\!=\{ (g,v,y)\in G\times\VV\times\YY\!\mid\! \mu(g,v)=y \}$. It fits into the diagram 
\[
\xymatrix@R=0.99em@C=3em{
&Z\ar[dl]^-\vpi\ar[dr]_-\psi&\\ G&&\VV\times\YY.}
\] 
If $(g,v,y)$ belongs to $Z$, then $v$, $y$ are in the same $G$-orbit and $\psi^{-1}(v,y)=G_{v,y}$. Let $Z_O\subset Z$ be the stratum corresponding to $v, y\in O$. Then $Z_O\srel{\psi}{\to}\VV_O\times\YY_O$ is smooth; its fibres are irreducible, isomorphic to $H$. Thus $Z_O$ is irreducible, and $\dim Z_O=\dim H+\dim\VV+\dim Y$. 

Now let $O'\cong G/H'$ be a lower-dimensional orbit in $\XX$; set $Z_{O'}:=\psi^{-1}(\VV_{O'}\times\YY_{O'})$. The fibres of $\psi|_{Z_{O'}}$ are isomorphic to $H'$. The assumptions on $\partial\VV$, $\partial\YY$ imply 
\begin{equation}
\begin{array}{r}
\dim Z_O-\dim Z_{O'} 
=\dim H'-\dim H\ges1.
\end{array}
\label{eq:OO}
\end{equation}
Hence $Z_O$ is the unique maximal dimensional, irreducible stratum of $Z$. To conclude that $Z$ is irreducible, we must show that $Z_{O'}$ is in the closure of $Z_O$. This follows again from~(\ref{iti:VY}): since $\psi|_{Z_{O'}}$ is an $H'$-fibre bundle, any $(g,v,y)\in Z_{O'}$ is the specialization of some $(g_1,v_1,y_1)\in\psi^{-1}(\VV_{O'}^{\rm smooth}\times\YY_{O'}^{\rm smooth})$. The transversality at~$v_1$ of the intersection $\VV\cap O'$ implies that $G\times\VV\srel{\mu}{\to}\XX$ is open about $(g_1,v_1)$. The (analytic) implicit function theorem shows that $(g_1,v_1,y_1)$ is the specialization of some $(g_2,v_2,y_2)\in Z_O$.

The morphism $Z\srel{\vpi}{\to} G$ is smooth at triples $(g,v,y)\in Z_O$ such that $g_*(\eT_{\VV,v})+\eT_{\YY,y}=\eT_{\XX,y};$ the corresponding fibre dimension is $\dim\VV+\dim\YY-\dim\XX>0$. So $\vpi(Z)\subset G$ is both closed ($\vpi$ is proper) and maximal dimensional; hence $\vpi$ is surjective. It follows that the intersection $g\VV\cap\YY=\vpi^{-1}(g)$ is non-empty for all $g\in G$.\smallskip

\emph{Claim}~2.~ The fibres of $\vpi$ are connected; for $g\in G^\gen$, $\vpi^{-1}(g)=g\VV\cap\YY$ is irreducible.

Let $Z'\srel{\nu}{\to} Z$ be the normalization; note that $Z'$ is irreducible. We consider the Stein factorization of $\vpi\nu$, $Z'\to G'\srel{\phi}{\to} G$. We prove that $\phi$ is \'etale; then $G'=G$, since $G$ is simply connected. Otherwise, the ramification locus of $\phi$ contains an irreducible divisor $D_{G'}\subset G'$; its pre-image in $Z'$ contains an irreducible, dominating divisor $D'$. 

\emph{Suppose $D'\cap Z_O\neq\emptyset$.}~ 
Since $Z_O\to\VV_O\times\YY_O$ is equidimensional, $D:=\overline{\psi(D')}\subset\VV\times\YY$ is either a divisor or the product. In both cases, $D\to\YY$ and $D\to\VV$ are surjective. Otherwise, one has $D=\VV\times D_\YY$, with $D_\YY\subset\YY$ an irreducible divisor (similarly for $D\to\VV$). By Claim~1, $D_\YY$ intersects every $g\VV$, so $D'\to G$ is surjective; this contradicts that $\vpi\nu(D')=\phi(D_{G'})$ is a divisor in $G$.
Consequently, for $(g,v,y)\in(\nu(D')\cap Z_O)^\gen$, we have $\codim_{G_{v,y}}(\nu(D')\cap G_{v,y})\les1$ and $\vpi$ is non-smooth along $\nu(D')\cap G_{v,y}$; this contradicts Inequality~\eqref{eq:p1}. Thus $\vpi\nu$ has connected fibres. For $g\in G^\gen$, ${(\vpi\nu)}^{-1}(g)$ is also normal; hence it is irreducible. The same holds for $\vpi$.

\emph{Suppose $D'\cap Z_O=\emptyset$.}~ There is a one-codimensional orbit $\BD$ such that $D'$ is (the closure of) a component of $\psi^{-1}(\VV_\BD\times\YY_\BD)$. We proved in Claim~1 that any translate of a component of $\YY_\BD$ intersects $\VV_\BD$. This contradicts Inequality~\eqref{eq:p1} again: for general $(g,v,y)\in\nu(D')\cap Z_\BD$, the codimension of the non-transverse intersection locus is at most one. 
\end{proof}

By contemplating Corollary~\ref{cor:barV} and Theorem~\ref{thm:mum}, one wonders if there is some transitivity property of the G3-property, that is, whether $\VV\cap\YY\subset\YY$ is still G3. The result below is in the same vein as~\cite{bad-deb}: it holds with more restrictive assumptions, but in a broader context; the results in \ocit\;concern arbitrary subvarieties of the product of two projective spaces only.

\begin{m-corollary}
Let $\VV$, $\YY$ be smooth with $2\codim_\XX\VV\les\dim\YY$, satisfying the hypotheses of Theorem~\ref{thm:mum}. In case~\eqref{iti:div}\eqref{itii:int}, we assume moreover that $[\YY]\cdot[\VV]^2\neq0$. Then, for $g\in G^\gen$, $g\VV\cap\YY$ is a smooth G3 $($actually $1^\pos)$ subvariety in $\YY$.
\end{m-corollary}

\begin{proof}
We consider the (possibly non-flat) movable family $W_g:=g\VV\cap\YY, g\in G$. Its general isotropy is the same as that of Diagram~\eqref{eq:strg-gen}, so it is irreducible. For $g\in G^\gen$, $W_g$ is connected, smooth (transverse intersection, by the $P(0)$-condition). Let $D_\YY\subset\YY$ be an irreducible divisor. 

In case~(\ref{iti:div})(\ref{itii:amp}), the normal bundle $\eN_{W_g/\YY}$ is a quotient of $\eN_{g\VV/\XX}$, thus ample. We deduce that $[W_g]^2\neq0$ and $W_g\subset\YY$ is G2. In  case~(\ref{iti:div})(\ref{itii:int}), $[W_g]^2\neq0$ is assumed, so the strong mobility (\textit{cf.} Lemma~\ref{lm:conn-mov}) yields the G2-property again. In both cases, we proved (see Claim~1) that $W_g\cap D_\YY\neq\emptyset$. The conclusion follows from Theorem~\ref{thm:g3}.
\end{proof}


\section{A criterion for complete intersections}

\begin{m-notation}\label{not:XV} 
Let $\YY$ be a smooth, irreducible subvariety of codimension two in a smooth projective variety~$\XX$, where $\dim\XX\ges4$, defined by the sheaf of ideals $\eI_\YY\subset\eO_\XX$. Its co-normal bundle $\eN_{\YY/\XX}^\vee:=\eI_\YY/\eI_\YY^2$ is locally free. Let $\kappa_X$ be the canonical line bundle of $\XX$.
\end{m-notation} 

\begin{m-question}
Suppose that 
$\eN_{\YY/\XX}^\vee=(\,\eL_1^{-1}\oplus\eL_2^{-1}\,)\otimes\eO_\YY$ for   
$\eL_1, \eL_2\in\Pic(\XX) $.
 
\nit Can one deduce that $\YY$ is a complete intersection? 
\end{m-question}
The interest stems from Hartshorne's conjecture, see~\cite{hart-conj}, that complete intersections are the only two-codimensional subvarieties of $\mbb P^n$, $n\ges6$. One may ask this for arbitrary ambient varieties. 

\begin{m-asmp}\label{asmp:yn}
Let the situation be as follows: 
\begin{enumerate}
\item\label{iti:llk} The variety $\XX$ is simply connected, and $\eL_1$, $\eL_2$, $\eL_2\eL_1^{-1}$ are semi-ample. 
\item\label{iti:mov} $\varstar$ 
$\YY$ is movable, belongs to a family $\cY$ with irreducible general isotropy, and intersects numerically non-trivially every effective divisor of $\XX$;

\nit $\varstar$ $\eL_1\rst_\YY$ is big and $\eL_2\rst_Y$ is ample. 
\item\label{iti:cod}
\emph{Either one} of the following conditions is satisfied: 
\begin{enumerate}
\item\label{itii:2cod} 
$\varstar$ $H^1(\XX,\eL_2\eL_1^{-1})=0$; 

\nit  $\varstar$ $\YY$ intersects non-trivially any two-codimensional subvariety of $X$.

\item\label{itii:fano} 
$\eL_1$, $\eL_2$ are ample, $\XX$ is Fano such that 
$\,H^t(\XX,\euf P)=0$ for all $\euf P\!\in\Pic(\XX)$, $t=1, 2, 3$.

\nit This is satisfied \textit{e.g.}~for $\XX$ with $\Pic(\XX)\cong\mbb Z$ and for products of such.
\end{enumerate}
\end{enumerate}
\end{m-asmp}

Let us comment on our assumptions. 
\begin{itemize}[leftmargin=3ex]
\item 
The subvariety $\YY$ is G3 (actually $1^\pos$) in $\XX$.
This follows from  Theorem~\ref{thm:g3}, because $\eN_{\YY/\XX}^\vee$ is not pseudo-effective, so $\YY$ is G2 in $\XX$. 

\item 
Let $\XX\srel{\vphi}{\to}\mbb P$ be the linear system of a power of either $\eL_1$ or $\eL_2$. Then its restriction $\vphi\rst_\YY$ is generically finite on image, and the Kodaira--Iitaka dimension of $\eL_1$ and $\eL_2$ is at least two.
\item[] 
Indeed, $\vphi^*\eO_{\mbb P}(1)$ is a power of $\eL_1$ (resp.~$\eL_2$). Since $\eL_1\rst_\YY$ is big, $\vphi\rst_\YY$ is generically finite on its image. (Otherwise, $\eL_1$, $\eL_2$ were trivial on movable curves in $\YY$.) Thus $\kod(\eL_1)$, $\kod(\eL_2)$ are bounded below by $\dim\vphi(\XX)\ges\dim\YY\ges2$. 

\item We analyse~(\ref{iti:cod})(\ref{itii:2cod}). 
\begin{itemize}[label=$\varstar$]
\item
The cohomology vanishing holds in the following (overlapping) situations:
\begin{itemize}[label=$\blacktriangleright$]
\item 
$\kappa_\XX^{-1}\eL_2\eL_1^{-1}$ is big, semi-ample: apply the Kawamata--Viehweg theorem; 
\item 
$\XX$ has the property $H^1(\XX,\euf P)=0$ for all nef $\euf P\in\Pic(\XX)$.  
\nit This is so for Fano varieties and many rational, almost homogeneous varieties ($G/P$, toric, etc.). Typically, they are Frobenius split, compatibly with an ample divisor. 
\end{itemize}

\item 
  As for intersections, there are at least two ways to ensure their non-triviality: 
\begin{itemize}[label=$\blacktriangleright$]
\item 
One is to impose $\YY$ to be a $2^\pos$ subvariety (\textit{cf.} \cite{hlc-pas}), so $\cd(\XX\sm\YY)\les\dim\XX-3$. Upper bounds for the cohomological dimension, in various situations, are given in~\cite{hlc+taj}. 
\item 
The second is suited when one knows the Chow/cohomology ring of $X$. Let $\Eff^{\mkern1mu 2}\!(\XX)$ be the convex cone generated by two-codimensional cycles of $\XX$ (in the 
vector space of numerical equivalence classes of algebraic cycles, with real coefficients); let $\overline{\Eff}^{\mkern1mu 2}\!(\XX)$ be its closure. For~(\ref{iti:cod})(\ref{itii:2cod}), it is enough to require that for all  
$c\in\Eff^{\mkern1mu 2}\!(\XX)\sm\{0\}$, we have $[\YY]\cdot c\neq0$.

A similar viewpoint is related to work of  Debarre--Ein--Lazarsfeld--Voisin~\cite{delv} and Fulger--Leh\-mann~\cite{ful+leh}. It suffices that $\YY$ contains a surface $M$ (\textit{e.g.}~complete intersection) whose class $[M]$ is in the interior of the cone  $\overline{\Eff}^{\mkern1mu 2}\!(\XX)^{\!\vee}$, the dual of the two-codimensional pseudo-effective cone. (We refer to \textit{op.\,cit.} for definitions.) 
\end{itemize}
\end{itemize}
\end{itemize}

\begin{m-theorem}\label{thm:g3-ci}
Let the situation be as in Assumptions~\ref{asmp:yn}. Then $\YY$ is a complete intersection in $\XX$.
\end{m-theorem}
Our approach is inspired from Faltings~\cite{falt-krit}, corresponding to the case of a projective space. 

\begin{proof}
Let $y_j\in\Gamma\big(\YY,(\eI_\YY/\eI_\YY^2)\otimes \eL_j\big)$, $j=1,2$, be the sections corresponding to the direct summand $\eO_\YY$. The proof consists in extending them to $\XX$. This is done in two steps: first we extend from $Y$ to $\hat\XX_\YY$, which requires the positivity assumptions on $\eL_1,\eL_2$; second, we further extend from $\hat\XX_\YY$ to $X$, which requires the G3-property.

\emph{Extension of $y_1$}.~ The homomorphisms 
$
\,\Gamma\big(\XX,(\eO_\XX/\eI_{\YY}^{r+1})\otimes \eL_{1}\big)\to 
\Gamma\big(\XX,(\eO_\XX/\eI_{\YY}^{r})\otimes \eL_{1}\big),\; r\ges 2,
$ 
are surjective (in fact, isomorphisms), so $y_1$ lifts to $\hat y_1\in\Gamma(\hat\XX_\YY, \hat\eL_1)$, where $\hat\eL_1\!:=\eL_{1}\otimes\eO_{\hat\XX_\YY}$. Indeed, we have the exact sequence 
\begin{m-eqn}{
0\lra \Sym^r\left(\frac{\eI_{\YY}}{\eI^2_{\YY}}\right)\otimes \eL_{1}
\lra\frac{\eO_\XX}{\eI^{r+1}_{\YY}}\otimes \eL_{1}
\lra \frac{\eO_\XX}{\eI^{r}_{\YY}}\otimes \eL_{1}\lra 0,
}\label{eq:NL}
\end{m-eqn}
where the first non-zero term is a direct sum of line bundles $(\euf M\otimes\eO_\YY)^{-1}$ of the form 
$$
\,\euf M=\eL_{1}^a\otimes\eL_{2}^b\otimes \eL_{1}^{-1},\quad a+b=r\ges2.
$$
These line bundles are big, so the Kawamata--Viehweg theorem yields the vanishing of the cohomology groups in degrees zero and one. 

We claim that $\Gamma(\XX,\eL_{1})\to \Gamma(\hat\XX_\YY, \hat\eL_1)$ is an isomorphism. Hartshorne~\cite[Proposition V.2.1]{hart-as} proved that G3-subvarieties possess the Lefschetz property: there is an open subset $U_\XX\subset\XX$ containing $\YY$, such that 
$\Gamma(U_\XX,\eL_{1}\rst_{U_\XX})\to \Gamma(\hat\XX_\YY,\hat\eL_1)$ 
is isomorphism. Let $z_{1,U_\XX}$ be the extension of $\hat y_1$; since $\XX\sm U_\XX\subset \XX\sm\YY$ contains no divisors, it extends to $z_1\in \Gamma(\XX,\eL_1)$. 

We obtained an element $z_1\in\Gamma(\XX, \eL_1)$ which restricts to $y_1\in\Gamma(\XX,\eI_\YY/\eI_\YY^2\otimes\eL_1)$. Its zero locus $Z$ is a divisor of $\XX$ containing $\YY$, and $y_1$ is actually the differential of $z_1$ along $\YY$. Since $\YY$ is smooth and $\eI_\YY/\eI_\YY^2$ splits, the differential criterion implies that $Z$ is smooth along $\YY$; hence there is only one irreducible component of $Z$ intersecting (containing) $\YY$. On the other hand, $\YY$ intersects the divisors of $\XX$, so all components of $Z$. Hence $Z$ is irreducible.

\emph{Extension of $y_2$}.~ We divide the discussion according to Assumptions~\ref{asmp:yn}(\ref{iti:cod}). 

\emph{Case}~(\ref{iti:cod})(\ref{itii:2cod}).~ The co-normal bundle of $\YY\subset Z$ is $(\eL_2\rst_\YY)^{-1}$, and exact sequences analogous to the sequence~\eqref{eq:NL} imply that $y_2$ lifts to $\hat y_2\in\Gamma(\hat Z_\YY,\eL_2)$. 
Note that $Y\subset Z$ is G3, being a divisor with ample normal bundle (\textit{cf.} \cite[Corollary 9.27]{bad}), so the Lefschetz property yields an open subset $U_Z\subset Z$ such that $\hat y_2$ extends to $y_{2,U_Z}'\in\Gamma(U_Z,\eL_2)$. But $\YY$ intersects all two-codimensional subvarieties of $\XX$, so $Z\sm U_Z$ contains no divisors of $Z$, and $y_{2,U_Z}'$ extends to $y_2'\in\Gamma(Z,\eL_2)$. The same argument as for $y_1$ shows that the zero locus of $y_2'$ is precisely $\YY$ (it has no further components). 
Finally, $y_2'$ extends to a section $z_2$ of $\eL_2$ over~$\XX$  because $H^1(\XX,\eL_2\eL_1^{-1})=0$; thus the restriction $\Gamma(X,\eL_2)\to\Gamma(Z,\eL_2)$ is surjective. 

\emph{Case}~(\ref{iti:cod})(\ref{itii:fano}).~ We repeat for $y_2$ the argument used for $y_1$ and obtain the sequence~\eqref{eq:NL} with $\eL_2$ instead of~$\eL_1$. To construct $\hat y_2\in\Gamma(\hat\XX_\YY,\hat\eL_2)$, it is necessary to have $H^1(\YY, \euf M^{-1})=0$ for 
\[
\euf M=\eL_{1}^a\otimes\eL_{2}^b\otimes\eL_{2}^{-1},\quad a+b\ges2.
\]
Once this is known, we extend $\hat y_2$ to a neighbourhood of $\YY$, then further to $z_2\in\Gamma(\XX,\eL_2)$. 

Let us prove the vanishing of the $H^1$-group. Since $\eL_1$ is ample, $Z$ has the property that $\Pic(\XX)\to\Pic(Z)$ is an isomorphism (\textit{cf.} \cite[Theorem IV.3.1]{hart-as}). From the sequence $0\to\eL_1^{-1}\to\eO_\XX\to\eO_Z\to0$, we deduce $H^t(Z,\euf P)=0$ for all $\euf P\in\Pic(Z)$, $t=1,2$. Finally, $Z$ is smooth along $\YY$, so $\eO_Z(-\YY)$  is locally free; we tensor $0\to\eO_Z(-\YY)\to\eO_Z\to\eO_\YY\to0$ by $\euf M^{-1}$ and obtain $H^1(\YY,\euf M^{-1})=0$. 
(Note: to conclude, we needed $H^1(Y,\euf M^{-1})=0$. One may ensure this by different assumptions, too, \textit{e.g.}~$(\eL_1^2\eL_2^{-1})\rst_\YY$ is semi-ample, big.)
\end{proof}

Conditions~(\ref{iti:cod})(\ref{itii:2cod}) are particularly suited for rational homogeneous varieties $\XX=G/P$, where $G$ is a semi-simple linear group and $P$ is parabolic. We assume $\ell_G\ges6$, where  
\[
\ell_G:=\,\text{minimum of the ranks of the simple factors of}\;G.
\]
First, the Barth--Lefschetz theorem, see~\cite{sena,soms-vdv}, implies that  $\Pic(\XX)\xrightarrow{\res}\Pic(\YY)$ is isomorphism, because $(\ell_G-2)/2\ges2$. Thus the summands of the split normal bundle $\eN_{\YY/\XX}$ are restricted from $\XX$. Second, Faltings~\cite[Satz 7]{falt-homog} proved that $\cd(\XX\sm\YY)\les\dim\XX-\ell_G+3<\dim\XX-2$; thus $\YY$ intersects every two-codimensional subvariety of $\XX$. (This still holds for $\YY$ lci. If $\eN_{\YY/\XX}$ is ample, the statement follows from the excess intersection theorem, see~\cite{fu-laz}, too.)  

\begin{m-corollary}\label{cor:rtl}
Let $\XX=G/P$ be as above, $\ell_G\ges6$, and let $\YY$ be a smooth, two-codimensional subvariety with split normal bundle $\eN_{\YY/\XX}=(\eL_1\oplus\eL_2)\otimes\eO_\YY$, satisfying  
\[
\eL_1,\;\eL_2,\;\eL_2\eL_1^{-1}\text{ semi-ample},\quad\eL_1\rst_\YY\text{ big},\quad\eL_2\rst_\YY\text{ ample.}
\]
Then $\YY$ is a complete intersection. If\, $G$ is simple and $P$ is maximal parabolic $($so $\Pic(\XX)\cong\mbb Z)$, the conditions on $\eL_1$, $\eL_2$ are automatically satisfied. 
\end{m-corollary}

\begin{proof} 
The higher cohomology of $\eL_2\eL_1^{-1}$ vanishes, because $\XX$ is homogeneous.
\end{proof}

\begin{m-corollary}\label{cor:fano}
Suppose $\XX$ is Fano, $\dim\XX\ges4$, with $\Pic(\XX)\cong\mbb Z$. Let $\YY\subset\XX$ be as in Assumptions~\ref{asmp:yn}\,\eqref{iti:llk}\eqref{iti:mov}, with $\eN_{\YY/\XX}^\vee=\big( \eO_\XX(-n_1)\oplus\eO_\XX(-n_2) \big)\rst_\YY$. Then $\YY$ is complete intersection. 
\end{m-corollary}

\begin{proof}
Here we are in situation~(\ref{iti:cod})(\ref{itii:fano}).
\end{proof}


\subsection{The lci case}\label{ssct:lci} 

The G2-criterion, see~\cite[Corollary 2.6]{hlc-g2}, and Theorem~\ref{thm:g3} are valid for lci subvarieties. Thus it is natural to ask what happens when $\YY\subset\XX$ is lci rather than smooth. 

The proof of Theorem~\ref{thm:g3-ci} shows that smoothness is involved at the following places: 
\begin{itemize}
\item 
The Kawamata--Viehweg theorem still yields the necessary cohomology vanishings, as long as the codimension of the singular locus of $\YY$ is at least two. (Note that lci implies Gorenstein, so Serre duality for $\YY$ works as in the smooth case.)
\item 
  The differential criterion for $z_1$ (resp.\ $y'_2$, $z_2$) extending $y_1$ (resp.\ $y_2$) still yields the smoothness of their zero loci along the regular locus $\YY\sm\Sing(\YY)$.

\nit  \emph{Claim}: If $\codim_\YY{\Sing}(\YY)\ges2$, the zero loci of $z_1$, $y'_2$, $z_2$ are irreducible. 

\nit  We consider $z_1$, which determines $Z$: let $Z_Y$ be the component containing $\YY$ and $Z'$ another component. Then $\YY\cap Z'\subset{\Sing}(\YY)$, so  
$2\les\codim_\YY{\Sing}(\YY)\les\codim_\YY(\YY\cap Z')\les1$. The other cases are similar.
\item 
The sheaf of ideals $\eI_{\YY\subset Z}=\eO_Z(-\YY)$ is locally free. This prevents us from considering  condition~(\ref{iti:cod})(\ref{itii:fano}) in this setting.
\item 
For homogeneous varieties $\XX$, the Barth--Lefschetz criterion for $\YY$ is not  available anymore. One has to impose that $\eL_1$, $\eL_2$ come from $\XX$, as in Assumption~\ref{asmp:yn}(\ref{iti:llk}). 
\end{itemize}

\begin{m-theorem}\label{thm:ci-sing}
Let Assumptions~\ref{asmp:yn}, with condition~\eqref{iti:cod}\eqref{itii:2cod}, be satisfied, and suppose that $\codim_\YY{\Sing(\YY)}\!\ges2$. Then Theorem~\ref{thm:g3-ci} remains valid; so does Corollary~\ref{cor:rtl}.
\end{m-theorem}


\subsection{A remark on Voisin's question}

The very moving property was considered in~\cite{voisin-coniv} to have geometric conditions which yield the bigness of a subvariety. 

\begin{conj-nono}[\textit{cf.} {\cite[Conjecture 2.5]{voisin-coniv}}]
If\, $Y\subset X$ is very moving $($possibly with ample normal bundle$)$, then $[Y]\in A_*(X)$ is big.
\end{conj-nono}

The conjecture holds for curves, divisors. A two-codimensional example shows that the ampleness of the normal bundle does not imply bigness. Theorem~\ref{thm:g3-ci} is a modest partial answer.

\begin{m-corollary}\label{cor:voisin}
Suppose Assumptions~\ref{asmp:yn} are satisfied and also that there exists an $n>0$ such that $\eL_1^n\eL_2^{-1}$ is semi-ample $(\eL_2$ is ``sandwiched'' between $\eL_1$, $\eL_1^n$, so $\eN_{\YY/\XX}$ is ample$)$. Then $\YY$ is big in $\XX$.
\end{m-corollary}

\begin{proof}
We proved that $\YY$ is the intersection of $Z_1\in|\eL_1|$ and  $Z_2\in|\eL_2|$. To deduce that $[\YY]$ is big, it is enough to prove that $\eL_1$, $\eL_2$ are so.
Let $\XX\xrightarrow{\vphi_1}\mbb P_1$ (resp.\ $\XX\xrightarrow{\vphi_2}\mbb P_2$) be the linear systems of suitable powers $\eL_1^a$ (resp.\ $\eL_2^b$). 

Note that $\YY$ is the pre-image by $Z_2\xrightarrow{\vphi_1\rst_{Z_2}}\mbb P_1$ of a hyperplane, and the pull-back of $\eO_{\mbb P_1}(1)$ is $\eL_1^a\rst_\YY$, which is big. So $\vphi_1\rst_\YY$ is finite on the image (otherwise, $\eL_1\rst_\YY$ is trivial on a movable curve), and the same holds for $\vphi_1\rst_{Z_2}$ (semi-continuity of fibre dimension). Hence $\eL_1\rst_{Z_2}$ is big. Similarly, $\eL_2\rst_{Z_1}$ is big. Now we use that $\eL_2=\eL_1\otimes(\eL_2\eL_1^{-1})$, where $\eL_2\eL_1^{-1}$ is semi-ample; therefore, $\eL_2\rst_{Z_2}$ is big, too. But $Z_2$ is the (reduced) pre-image of a hyperplane in $\mbb P_2$, so $\vphi_2\rst_{Z_2}$ is generically finite on the image. Therefore, $\vphi_2$ has the same property, so $\eL_2$ is big. 

We use the ``sandwich'' property for dealing with $\eL_1$. The bigness of $\eL_2\rst_{Z_1}$ yields that of $\eL_1\rst_{Z_1}$, which in turn implies that $\vphi_1$ generically finite. It follows that $\eL_1$ is big. 
\end{proof}


\newcommand{\etalchar}[1]{$^{#1}$}

\end{document}